\documentclass[12pt,reqno,a4paper]{amsart}
\usepackage{placeins}
\usepackage{colonequals}
\usepackage{booktabs}
\usepackage{colonequals}
\usepackage[T1]{fontenc}
\usepackage[utf8]{inputenc}
\usepackage{amsmath,amssymb,amsthm}
\usepackage{enumitem}
\usepackage{a4wide}
\usepackage{tikz}
\usetikzlibrary{positioning}
\usepackage[colorlinks, hyperindex, plainpages=false, bookmarksopen, bookmarksnumbered, pdfusetitle,  linkcolor=blue, citecolor=blue ]{hyperref}
\setlist{itemsep=1pt,topsep=3pt,parsep=0pt}
\usepackage[a4paper,hmargin=1.9cm,top=1.6cm,bottom=1.7cm]{geometry}
\theoremstyle{plain}
\newtheorem{theorem}{Theorem}[section]
\newtheorem{lemma}[theorem]{Lemma}
\newtheorem{proposition}[theorem]{Proposition}
\newtheorem{corollary}[theorem]{Corollary}
\newtheorem{mainthm}{Theorem}

\theoremstyle{definition}
\newtheorem{definition}[theorem]{Definition}
\newtheorem{example}[theorem]{Example}
\newtheorem{remark}[theorem]{Remark}

\newcommand{\CC}{\mathbb{C}}
\newcommand{\ZZ}{\mathbb{Z}}
\newcommand{\TT}{\mathbb{T}}
\newcommand{\DD}{\mathbb{D}}
\newcommand{\RR}{\mathbb{R}}

\newcommand{\J}{\mathbf{J}}
\newcommand{\E}{\mathcal{E}}
\newcommand{\diag}{\mathrm{diag}}
\newcommand{\WLip}{\mathcal W^{\mathrm{Lip}}_{\mathrm{Dou}}}
\newcommand{\WLtwo}{\mathcal W^{L^2}_{\mathrm{Dou}}}
\newcommand{\WDou}{\mathcal W_{\mathrm{Dou}}}
\DeclareMathOperator{\supp}{supp}
\renewcommand{\Re}{\mathrm{Re}\,}
\renewcommand{\Im}{\mathrm{Im}\,}

\title{A Characterization of Douglas Weights}

\author{Y.~Elmadani}
\author{A.~Hanine}
\address{LAMA, Mohammed V University in Rabat, Faculty of Sciences, Rabat Morocco}
\email{elmadanima@gmail.com}
\email{abhanine@gmail.com}

\begin{document}
	
			\maketitle
		
	\begin{abstract}
		A weight $\omega$ on the unit disc is a \emph{Douglas weight} if the $\omega$-energy of the Poisson
		extension of every Lipschitz function $f$ on the unit circle equals
		$\iint|f(\zeta)-f(\lambda)|^2\,d\J_\omega$ for some positive symmetric measure $\J_\omega$ off
		the diagonal. For radial weights, we show that this holds if and only if	\[
		j_\omega(\theta)
		=
		\int_0^1 \frac{\omega(r)}{r}\,
		\partial_\theta^2 P(r^2,\theta)\,dr
		\geq 0,
		\qquad 0<\theta\leq\pi,
		\]
		where $P(x,\theta)$ denotes the Poisson kernel. We also obtain a complete characterization of the Douglas weights among arbitrary weights. The class of Douglas weights
		strictly contains the nonnegative superharmonic weights, and we give examples that are
		neither superharmonic nor subharmonic.
	\end{abstract}
	\maketitle
	\tableofcontents
	\section{Introduction}\label{sec:intro}
	
	Let $\DD$ be the unit disc and $\TT$ the unit circle, with the normalized area measure $dA=\frac1\pi\,dx\,dy$ on $\DD$ and the normalized Lebesgue measure $m$ on $\TT$. For a holomorphic function $F$ in $\DD$ with boundary function $f$, the Dirichlet integral $D(F)=\int_{\DD}|F'|^2\,dA$ is a quadratic form in $f$, with an explicit kernel acting on the differences $f(\zeta)-f(\lambda)$:
	\begin{equation}\label{eq:intro-douglas}
		D(F)=\iint_{\TT\times\TT}\frac{|f(\zeta)-f(\lambda)|^2}{|\zeta-\lambda|^2}\,dm(\zeta)\,dm(\lambda)=\sum_{n\in\ZZ}|n|\,|\hat f(n)|^2 .
	\end{equation}
	This identity is due to Douglas \cite{douglas1931solution}; see also
	\cite{carleson1960representation,el2014primer,koosis1998introduction}.
	
	For a weight $\omega\ge 0$ on $\DD$, the weighted Dirichlet integral of a
	holomorphic function $F$ on $\DD$ is
	$D_\omega(F)=\int_{\DD}|F'|^2\,\omega\,dA$. When $\omega=P[\nu]$ is the Poisson integral of a finite positive
	measure $\nu$ on $\TT$, Richter and Sundberg \cite{richter1991formula} proved a
	Douglas-type formula, which yields
	\[
	D_\omega(F)=\int_{\TT} D_\zeta(F)\,d\nu(\zeta),
	\qquad
	D_\zeta(F)=\int_{\TT}\frac{|f(\lambda)-f(\zeta)|^2}{|\lambda-\zeta|^2}\,dm(\lambda),
	\]
	where $f(\lambda)$ denotes the nontangential boundary values of $F$ and $D_\zeta(F)=\infty$ if $F$ has no
	nontangential limit at $\zeta$. The associated spaces $\mathcal D_\omega$ were introduced by Richter
	\cite{richter1991representation}, who showed that every cyclic analytic
	two-isometry is unitarily equivalent to multiplication by $z$ on some
	$\mathcal D_{P[\nu]}$. Aleman \cite{aleman1993multiplication} extended much
	of this theory to positive superharmonic weights. Bao, G\"o\u{g}\"u\c{s} and Pouliasis \cite{bao2018dirichlet} proved a Douglas-type formula for superharmonic weights that are Green potentials of measures on $\DD$. For weighted Besov spaces, Bahajji-El~Idrissi and El-Fallah
	\cite{idrissi2025douglas} obtained a two-sided Douglas-type estimate,
	valid up to multiplicative constants. Dirichlet spaces with harmonic or superharmonic weights also appear in connection with de Branges--Rovnyak spaces \cite{Sarason1997,ChevrotGuillotRansford2010,el2016dirichlet}, reproducing kernels and the complete Nevanlinna--Pick property \cite{shimorin2001reproducing,shimorin2002complete}, invariant subspaces and cyclicity \cite{ross1994hyperinvariant,el2009brown,el2016cyclicity,elmadani2019cyclicity,abakumov2015cyclicity,el2022extremal}, boundary behaviour \cite{guillot2012fine,el2024ahlfors}, zero sets, approximation and capacity \cite{idrissi2020blaschke,bahajji2020approximation,el2019kernel}, and zero-free regions of the Riemann zeta function \cite{gallardo2025zero}.
	
	Formulas of this type also arise in the theory of Dirichlet forms. Let $\J$ be a positive symmetric measure on pairs of distinct points of $\TT$ with $\iint|\zeta-\lambda|^2\,d\J<\infty$, and consider the form $f\mapsto\iint|f(\zeta)-f(\lambda)|^2\,d\J$ on Lipschitz functions. When this form is closable in $L^2(\TT)$, its closure is a regular Dirichlet form of pure jump type in the sense of Beurling and Deny \cite{beurling1958espaces,fukushima2010dirichlet}, associated with a symmetric pure-jump Hunt process on the circle; the kernel $|\zeta-\lambda|^{-2}$ of \eqref{eq:intro-douglas} gives the Cauchy process \cite{bertoin1996levy}. The companion paper \cite{BEHJL2026} studies Douglas weighted Dirichlet spaces from this point of view. The trace on $\TT$ of the weighted energy $u\mapsto\int_\DD|\nabla u|^2\omega\,dA$, studied in \cite{fukushima2004timechanges,chen2006traces,chenfukushima2012}, minimizes the energy over all extensions of a boundary function and is always a Dirichlet form.
	
	Throughout the paper, a \emph{weight} is a measurable function $\omega\ge0$ on $\DD$ with $\|\omega\|_1=\int_\DD\omega\,dA<\infty$. For $f\in L^2(\TT)$ let $P[f]$ be its Poisson extension and
	$$\E_\omega(f)=\int_{\DD}|\nabla P[f]|^2\,\omega\,dA\in[0,\infty],\qquad|\nabla u|^2=|u_x|^2+|u_y|^2 .$$
	If $f$ is the boundary function of a holomorphic function $F$, then $|\nabla F|^2=2|F'|^2$, so $\E_\omega(f)=2D_\omega(F)$. Let $\diag=\{(\zeta,\zeta):\zeta\in\TT\}$. We say that $\omega$ is a \emph{Douglas weight} if there is a positive symmetric Borel measure $\J_\omega$ on $\TT^2\setminus\diag$ such that
	\begin{equation}\label{eq:intro-def}
		\E_\omega(f)=\iint_{\TT\times\TT}|f(\zeta)-f(\lambda)|^2\,d\J_\omega(\zeta,\lambda)
	\end{equation}
	for every Lipschitz function $f$ on $\TT$, and we write $\WDou=\WLip$ for the class of Douglas weights. If there is such a measure for which \eqref{eq:intro-def} holds for every Borel function $f$ on $\TT$ with $\int_\TT|f|^2\,dm<\infty$, we call $\omega$ a \emph{strong Douglas weight}, and we write $\WLtwo$ for the class of strong Douglas weights. Thus $\WLtwo\subset\WLip$, and both classes are convex cones, since $\E_\omega(f)$ is linear in $\omega$ and the right-hand side of \eqref{eq:intro-def} is linear in $\J_\omega$. The left-hand side of \eqref{eq:intro-def} depends only on the class of $f$ in $L^2(\TT)$, whereas the right-hand side may depend on the values of $f$ on an $m$-null set. The two classes coincide for radial weights (Theorem~\ref{thm:A}), but not in general: the Poisson kernel $P_{\zeta_0}$ is a Douglas weight that is not a strong Douglas weight; see Proposition~\ref{prop:strict}. Because \eqref{eq:intro-def} is required for all Lipschitz functions, the measure $\J_\omega$ is unique (see the end of the proof of Theorem~\ref{thm:A}); with boundary values of holomorphic functions only, uniqueness fails (Remark~\ref{rem:holomorphic}).
	
	Every positive superharmonic weight is a Douglas weight, with an explicit measure. Let $\omega\not\equiv+\infty$ be a positive superharmonic function on $\DD$. By the Riesz decomposition theorem (see \cite{aleman1993multiplication,el2016dirichlet,bahajji2026superharmonically}),
	$$\omega=G[\mu]+P[\nu],\qquad G[\mu](z)=\int_\DD2\log\Big|\frac{1-\bar wz}{z-w}\Big|\,d\mu(w),$$
	where $\mu=-\frac1{4\pi}\Delta\omega\,dx\,dy$ is a positive measure on $\DD$ with $\int_\DD(1-|w|^2)\,d\mu(w)<\infty$, the Riesz measure of $\omega$, and $\nu$ is a finite positive measure on $\TT$, the boundary measure of $\omega$. With this normalization $G[dA](z)=1-|z|^2$, and Tonelli's theorem together with the symmetry of the Green function gives $\int_\DD G[\mu]\,dA=\int_\DD(1-|w|^2)\,d\mu(w)$. Since $\int_\DD P[\nu]\,dA=\nu(\TT)$, the function $\omega$ is integrable. For $\zeta\ne\lambda$ in $\TT$ we put
	\begin{equation}\label{eq:Amu}
		A_\mu(\zeta,\lambda)=\int_\DD P_\zeta(z)\,P_\lambda(z)\,d\mu(z)\in[0,\infty].
	\end{equation}
	
	\begin{proposition}\label{prop:superharmonic}
		Let $\omega$ be a positive superharmonic function on $\DD$, with Riesz measure $\mu$ and boundary measure $\nu$. Then, for every $f\in\mathrm{Lip}(\TT)$,
		\begin{equation}\label{Douglas}
			\int_{\DD}|\nabla P[f]|^2\,\omega\,dA=\iint_{\TT\times\TT}|f(\zeta)-f(\lambda)|^2\left(A_\mu(\zeta,\lambda)\,dm(\zeta)\,dm(\lambda)+\frac{2\,d\nu(\zeta)\,dm(\lambda)}{|\zeta-\lambda|^2}\right).
		\end{equation}
		In particular $\omega$ is a Douglas weight, with
		$$\J_\omega=A_\mu\,dm\,dm+\frac{\nu\otimes m+m\otimes\nu}{|\zeta-\lambda|^2}.$$
	\end{proposition}
	
Proposition~\ref{prop:superharmonic} follows from
\cite{richter1991formula} for harmonic weights,
\cite{bao2018dirichlet} for Green potentials, and
\cite[Section~4]{bahajji2026superharmonically} for general
superharmonic weights. These works consider holomorphic
functions and use the normalization
$D_\omega(F)=\frac12\E_\omega(f)$. The cross terms vanish
for trigonometric polynomials, and
Lemma~\ref{lem:approx} extends the identity to Lipschitz
functions.
	
	The converse of Proposition~\ref{prop:superharmonic} fails: not every Douglas
	weight is superharmonic; see Example~\ref{prop:subharm}. Moreover, the class of
	Douglas weights is not solid, that is, a weight dominated by a Douglas weight
	need not be a Douglas weight. Indeed, a nonzero bounded weight $\omega_0$ with
	compact support in $\DD$ is never a Douglas weight, whereas $\omega_0+c$ is
	one for every sufficiently large constant $c$; see
	Proposition~\ref{prop:nonsolid}.
	
	The aim of this paper is to characterize Douglas weights.\\
	
	Let $\tilde\omega(r)=\frac1{2\pi}\int_0^{2\pi}\omega(re^{it})\,dt$ be the radialization of $\omega$, and let $P(x,\theta)=\frac{1-x^2}{1-2x\cos\theta+x^2}$ for $0\le x<1$. For $\theta\in\RR\setminus2\pi\ZZ$ put
	\begin{equation}\label{eq:jomega}
		j_\omega(\theta)=\int_0^1\frac{\tilde\omega(r)}{r}\,\partial_\theta^2P(r^2,\theta)\,dr .
	\end{equation}
	The integral converges absolutely; see \eqref{eq:Ctheta}. For the constant weight $\omega\equiv 1$, we have
	\[
	j_1(\theta)=\frac{1}{4\sin^2(\theta/2)},
	\]
	and hence
	\[
	j_1(\arg(\bar\zeta\lambda))=|\zeta-\lambda|^{-2},
	\]
	which recovers the kernel in \eqref{eq:intro-douglas}; see Example~\ref{ex:j1}.
	
	\begin{mainthm}\label{thm:A}
		Let $\omega\ge0$ be a radial weight, $\omega(z)=\omega(|z|)$. The following are equivalent.
		\begin{enumerate}
			\item[(i)] $\omega\in\WLtwo$, that is, $\omega$ is a strong Douglas weight.
			\item[(i')] $\omega\in\WLip$, that is, $\omega$ is a Douglas weight.
			\item[(ii)] $j_\omega(\theta)\ge0$ for every $\theta\in(0,\pi]$.
		\end{enumerate}
		In that case $\J_\omega(d\zeta,d\lambda)=2\,j_\omega(\arg\bar\zeta\lambda)\,dm(\zeta)\,dm(\lambda)$.
	\end{mainthm}
	
	For the superharmonic weight $\omega(z)=1-|z|^2$, the
	corresponding function $j_\omega$ admits the explicit
	expression
	\[
	j_{1-|z|^2}(\theta)
	=
	-\frac12
	-\cos\theta\,\log\left(2\sin\frac{\theta}{2}\right)
	+\frac12(\pi-\theta)\sin\theta,
	\qquad 0<\theta<2\pi.
	\]
	Moreover, for every $F\in H^2$ with boundary function $f$,
	we have
	\begin{equation}\label{eq:intro-H2}
		\int_{\mathbb D}|F'|^2(1-|z|^2)\,dA
		=
		\iint_{\mathbb T\times\mathbb T}
		|f(\zeta)-f(\lambda)|^2
		j_{1-|z|^2}(\arg(\bar\zeta\lambda))
		\,dm(\zeta)\,dm(\lambda).
	\end{equation}
	See Example~\ref{lem:AdA}.
	Formula \eqref{eq:intro-H2} should be compared with a result of Bahajji-El~Idrissi and El-Fallah \cite{idrissi2025douglas}. For $\varpi:(0,1]\to[0,\infty)$ and $1<p<\infty$, they consider $\int_\DD|F'|^p\varpi(1-|z|^2)\,dA$ and put $\varpi_p(x)=\varpi(x)/x^p$. They call $\varpi$ \emph{$p$-admissible} if $t\mapsto t^\beta\varpi(t)$ is non-decreasing for some $\beta\in(0,1)$ and $\int_x^1\varpi_p(t)\,dt\le Cx\,\varpi_p(x)$ for all $x\in(0,1)$, and they prove \cite[Theorem~1.1]{idrissi2025douglas} that, for such $\varpi$ and $F\in H^p$,
	$$\|F\|_{H^p}^p+\int_\DD|F'|^p\varpi(1-|z|^2)\,dA\asymp\|F\|_{H^p}^p+\iint_{\TT\times\TT}|f(\zeta)-f(\lambda)|^p\,\varpi_p(|\zeta-\lambda|)\,dm(\zeta)\,dm(\lambda).$$
	The weight $1-|z|^2$ corresponds to $p=2$ and $\varpi(t)=t$. Then $\varpi_2(t)=1/t$, and $\int_x^1\varpi_2(t)\,dt=\log(1/x)$ is not bounded by $Cx\,\varpi_2(x)=C$, so $\varpi$ is not $2$-admissible. This weight is the endpoint $\alpha=p-1$ of the power weights $t^\alpha$, which are $p$-admissible exactly for $-1<\alpha<p-1$, and \cite{idrissi2025douglas} shows with the functions $1/(1-az)$ that the two-sided estimate fails at this endpoint. In \eqref{eq:intro-H2} the kernel grows only logarithmically at the diagonal, since $j_{1-|z|^2}(\theta)=\log(1/|\theta|)-\frac12+o(1)$ as $\theta\to0$, while $\varpi_2(|\zeta-\lambda|)=|\zeta-\lambda|^{-1}$. Thus $1-|z|^2$ admits an exact Douglas formula, although it is not $2$-admissible. The estimate of \cite{idrissi2025douglas} holds for every $1<p<\infty$, for holomorphic functions, with the kernel $\varpi_p(|\zeta-\lambda|)$ determined by $\varpi$; Theorem~\ref{thm:A} concerns $p=2$ and decides, for every radial weight, whether an exact formula \eqref{eq:intro-def} exists, and gives its kernel.
	
	For general weights, our characterization is expressed through integrals of
	$\omega$ against products of two Poisson kernels. For $\zeta\in\TT$ let
	$P_\zeta(z)=(1-|z|^2)/|\zeta-z|^2$, and for $\varphi\in C(\TT^2)$ define
	\begin{equation}\label{eq:intro-Pphi}
		\mathcal P\varphi(z)=\iint_{\TT\times\TT}\varphi(\zeta,\lambda)\,
		P_\zeta(z)\,P_\lambda(z)\,dm(\zeta)\,dm(\lambda)\qquad(z\in\DD).
	\end{equation}
	If $\varphi\in C^\infty(\TT^2)$, then $\mathcal P\varphi$ is smooth in $\DD$
	and $\Delta\mathcal P\varphi$ is bounded on $\DD$; see Lemma~\ref{lem:elem}.
	
	\begin{definition}\label{def:weaksuper}
		Let $\omega\ge0$ be a measurable function on $\DD$ such that $\int_\DD\omega\,dA<\infty$. We say that $\omega$ \emph{satisfies condition~\textup{(EH)}} if
		\begin{equation}\label{eq:weaksuper}
			\int_\DD\omega\,\Delta\big(\mathcal P\varphi\big)\,dA\le0
			\tag{EH}
		\end{equation}
		for every nonnegative function $\varphi\in C_c^\infty(\TT^2\setminus\diag)$. Since $\Delta(\mathcal P\varphi)$ is bounded on $\DD$, the integral in \eqref{eq:weaksuper} converges absolutely.
	\end{definition}
	
	\begin{mainthm}\label{thm:B}
		Let $\omega\ge0$ be a measurable function on $\DD$ such that $\int_\DD\omega\,dA<\infty$. Then $\omega$ is a Douglas weight if and only if $\omega$ satisfies condition~\eqref{eq:weaksuper}.
	\end{mainthm}
	
	The measure $\J_\omega$ is then given by
	$\iint\varphi\,d\J_\omega=-\tfrac14\int_\DD\omega\,\Delta(\mathcal P\varphi)\,dA$;
	see Theorem~\ref{thm:main}. Every positive superharmonic weight satisfies
	\eqref{eq:weaksuper}, but this condition is strictly weaker than
	$\Delta\omega\le0$. Indeed, for nonzero $\varphi\ge0$ the function
	$\mathcal P\varphi$ is positive throughout $\DD$, so \eqref{eq:weaksuper} is not
	a local condition: it also involves the behavior of $\omega$ near $\TT$. For
	instance, the weight $\omega(z)=|z|^2$ satisfies \eqref{eq:weaksuper},
	although $\Delta|z|^2=4$; see Remark~\ref{rem:weak}.
	
	In Sections~\ref{sec:further} and~\ref{sec:examples} we apply the two theorems. The weight $(1-|z|^2)^s$ is a Douglas weight exactly when $-1<s\le1$, and every integrable mixture $\int|z|^{2\gamma}\rho(d\gamma)$ with $\rho\ge0$ is a Douglas weight, see Examples~\ref{prop:subharm} and~\ref{prop:power}. A radial weight $1+\varphi$ with $\|\varphi\|_\infty\le\frac45$ is a Douglas weight; see Theorem~\ref{thm:ball}, and this gives smooth radial Douglas weights that are neither superharmonic nor subharmonic; see Example~\ref{ex:wild}. No such bound exists for non-radial perturbations: for every $\varepsilon\in(0,1]$ there is a weight $1+\varepsilon\phi$ with $|\phi|\le1$ that is not a Douglas weight; see Proposition~\ref{prop:nostab}. For a bounded $C^2$ weight with bounded Laplacian, Theorem~\ref{thm:B} becomes a pointwise inequality between $\Delta\omega$ and the boundary values of $\omega$; for instance, $1+\varepsilon\,r^3\cos t$ is a Douglas weight for $0<\varepsilon\le\frac13$; see Example~\ref{ex:r3cos}. The weights $1+\varepsilon\cos(\arg z)$, $0<\varepsilon\le1$, have constant radialization; they are Douglas weights exactly when $\varepsilon\le\varepsilon^*$, for a number $\varepsilon^*<1$ whose numerical value is about $0.834$ (Example~\ref{ex:cos}). Finally, the Poisson kernel $P_\zeta$ is a Douglas weight for which \eqref{eq:intro-def} cannot hold for all Borel functions $f\in L^2(\TT)$; see Proposition~\ref{prop:strict}.
	
The paper is organized as follows. Section~\ref{sec:radial} contains the proof of Theorem~\ref{thm:A}, together with the Fourier analysis of the Poisson kernel that it requires. Section~\ref{sec:main} contains the proof of Theorem~\ref{thm:B}, based on the matrix of the energy form in the Fourier basis and on the distribution $K_\omega$ on $\TT^2$ defined by the energy form; it uses some identities established in Section~\ref{sec:radial}. Section~\ref{sec:further} contains further results: necessary conditions, negative definite functions, weights with a bounded Laplacian and Dirichlet forms. Section~\ref{sec:examples} is devoted to examples and counterexamples.\\
	
	\paragraph{\textit{Notation and conventions}}\label{subsec:setting}
	We write $e_n(\zeta)=\zeta^n$ for $n\in\ZZ$, identify $\TT$ with $\RR/2\pi\ZZ$ through $\theta\mapsto e^{i\theta}$, and write $\TT\setminus\{0\}$ for the circle without the point $\theta=0$. For a measurable function $F\ge0$ on $\DD$,
	\begin{equation}\label{eq:polar}
		\int_{\DD}F\,dA=2\int_0^1 r\Big(\frac1{2\pi}\int_0^{2\pi}F(re^{it})\,dt\Big)dr .
	\end{equation}
	The angular modes of a weight are $\omega_j(r)=\frac1{2\pi}\int_0^{2\pi}\omega(re^{it})e^{-ijt}\,dt$, so that $\omega_0=\tilde\omega$, $\omega_{-j}=\overline{\omega_j}$ and $|\omega_j|\le\tilde\omega$; by \eqref{eq:polar}, $\int_0^1\tilde\omega(r)\,r\,dr=\frac12\|\omega\|_1$. For a function $f$ and a measure $\nu$ on $\TT$, $P[f]$ and $P[\nu]$ denote their Poisson integrals.
	
	For $f\in L^2(\TT)$ we write $Q_\omega(f)=\E_\omega(f)$, and we use the sesquilinear form $\E_\omega(f,g)=\int_\DD\nabla P[f]\cdot\overline{\nabla P[g]}\,\omega\,dA$ and the bilinear form $Q_\omega(g,h)=\int_\DD\nabla P[g]\cdot\nabla P[h]\,\omega\,dA$, where $\nabla u\cdot\nabla v=u_xv_x+u_yv_y$. Thus $\E_\omega(f,g)=Q_\omega(f,\bar g)$ and $Q_\omega(f)=Q_\omega(f,\bar f)$. The matrix $M_\omega(n,k)=\E_\omega(e_n,e_k)=Q_\omega(e_n,e_{-k})$ is computed in Lemma~\ref{lem:matrix}. A positive symmetric Borel measure $\J$ on $\TT^2\setminus\diag$ defines, for Borel functions on $\TT$,
	$$Q_\J(f)=\iint|f(\zeta)-f(\lambda)|^2\,d\J\in[0,\infty],\qquad Q_\J(g,h)=\iint\big(g(\zeta)-g(\lambda)\big)\big(h(\zeta)-h(\lambda)\big)\,d\J,$$
	the second one when $Q_\J(g)$ and $Q_\J(h)$ are finite. Finally, $\mathrm{Lip}(f)=\sup_{\zeta\ne\lambda}|f(\zeta)-f(\lambda)|/|\zeta-\lambda|$, $D^h(\omega)=\{f\in L^2(\TT):Q_\omega(f)<\infty\}$, $D(\J)=\{f\in L^2(\TT):Q_\J(f)<\infty\}$, and $F_\omega$ is the closure of $C^\infty(\TT)$ in the norm $(\|f\|_{L^2}^2+Q_\omega(f))^{1/2}$.
	
	\section{Proof of Theorem~\ref{thm:A}}\label{sec:radial}
	
	We prove the implications $(ii)\Rightarrow(i)$,
	$(i)\Rightarrow(i')$, and $(i')\Rightarrow(ii)$ in
	Theorem~\ref{thm:A}, and obtain an explicit formula for
	$\J_\omega$. We also discuss the relation between the
	domains of the two forms. The proof is based on their
	Fourier representations. For radial weights, both forms
	are diagonal in the Fourier basis. The function
	$j_\omega$ is determined by the Fourier coefficients of
	the second angular derivative of the Poisson kernel.
	Several identities derived in the proof hold for arbitrary
	weights and will be used again in
	Sections~\ref{sec:main} and~\ref{sec:examples}.
	
	\begin{proof}[Proof of Theorem~\ref{thm:A}]
		We begin with the energy of the exponentials, for an arbitrary weight $\omega$. For $n\ge1$ the function $u=P[e_n]=z^n$ is holomorphic, so $u_x=nz^{n-1}$, $u_y=inz^{n-1}$ and $|\nabla u|^2=2n^2|z|^{2n-2}$; for $n\le-1$, $P[e_n]=\bar z^{|n|}$ and the computation is the same. By \eqref{eq:polar},
		\begin{equation}\label{eq:psi}
			\psi_\omega(n):=Q_\omega(e_n)=2n^2\int_\DD|z|^{2|n|-2}\omega\,dA=4n^2\int_0^1r^{2|n|-1}\,\tilde\omega(r)\,dr\quad(n\ne0),\qquad\psi_\omega(0)=0 .
		\end{equation}
		Thus $\psi_\omega$ depends on $\omega$ only through $\tilde\omega$, $\psi_\omega(\pm1)=2\|\omega\|_1$, and $\psi_\omega(n)\le2n^2\|\omega\|_1$ because $|z|^{2|n|-2}\le1$. When $\omega$ is radial, the energy of every $f\in L^2(\TT)$ is diagonal in the Fourier basis. Indeed, for $u=P[f]=\sum_n\hat f(n)r^{|n|}e^{int}$ we have $|\nabla u|^2=|u_r|^2+r^{-2}|u_t|^2$, Parseval's identity on the circle $\{|z|=r\}$ gives $\frac1{2\pi}\int_0^{2\pi}|\nabla u(re^{it})|^2dt=\sum_n2n^2r^{2|n|-2}|\hat f(n)|^2$, and \eqref{eq:polar} yields
		\begin{equation}\label{eq:parseval}
			Q_\omega(f)=\sum_{n\in\ZZ}\psi_\omega(n)\,|\hat f(n)|^2\in[0,\infty]\qquad(\omega\ \text{radial},\ f\in L^2(\TT)).
		\end{equation}
		
		To bring in $j_\omega$ we need the second angular derivative of the Poisson kernel. Write $P(x,\theta)=(1-x^2)/D$ for $0\le x<1$, where $D=D(x,\theta)=1-2x\cos\theta+x^2=\delta^2+2x(1-\cos\theta)$ and $\delta=1-x$, and put $q=xe^{i\theta}$. We use the classical expansion \cite{koosis1998introduction}
		\begin{equation}\label{eq:Pseries}
			P(x,\theta)=\sum_{k\in\ZZ}x^{|k|}e^{ik\theta},
		\end{equation}
		which gives $\int_{-\pi}^{\pi}P(x,\theta)\cos n\theta\,d\theta=2\pi x^{|n|}$ by orthogonality. Differentiating \eqref{eq:Pseries} twice in $\theta$, and differentiating $(1-x^2)/D$ twice, we obtain
		\begin{equation}\label{eq:d2Pseries}
			\partial_\theta^2P(x,\theta)=-\sum_{k\in\ZZ}k^2x^{|k|}e^{ik\theta}=-2\Re\frac{q(1+q)}{(1-q)^3}=\frac{2x(1-x^2)\,N}{D^3},\qquad N=4x\sin^2\theta-D\cos\theta .
		\end{equation}
		Since $\sin^2\theta=(1-\cos\theta)(1+\cos\theta)\le2(1-\cos\theta)$ and $D\ge2x(1-\cos\theta)$, we have $4x\sin^2\theta\le4D$, hence $|N|\le5D$; as $1-x^2=\delta(1+x)\le2\delta$, this gives $|\partial_\theta^2P|\le10\,x(1-x^2)D^{-2}\le20\,\delta D^{-2}$. If $x\le\frac12$, then $D\ge\delta^2\ge\frac14$. If $x\ge\frac12$ and $|\theta|\le\pi$, then $D\ge\delta^2+(1-\cos\theta)\ge\frac2{\pi^2}(\delta^2+\theta^2)$, because $1-\cos\theta=2\sin^2\frac\theta2\ge2\theta^2/\pi^2$. Therefore, for $|\theta|\le\pi$,
		\begin{equation}\label{eq:d2Pbound}
			|\partial_\theta^2P(x,\theta)|\le160\,x\quad\big(x\le\tfrac12\big),\qquad|\partial_\theta^2P(x,\theta)|\le\frac{5\pi^4\,\delta}{(\delta^2+\theta^2)^2}\quad\big(x\ge\tfrac12\big).
		\end{equation}
		Two consequences of \eqref{eq:d2Pbound} control $j_\omega$. First, for fixed $\theta\in[-\pi,\pi]\setminus\{0\}$ the second bound is at most $5\pi^4\theta^{-4}\le10\pi^4\theta^{-4}x$, so
		\begin{equation}\label{eq:Ctheta}
			|\partial_\theta^2P(x,\theta)|\le C(\theta)\,x\quad(0\le x<1),\qquad C(\theta)=\max(160,\,10\pi^4\theta^{-4}),
		\end{equation}
		and $C$ is bounded on compact subsets of $\TT\setminus\{0\}$. Second, for $n\ge1$, using $1-\cos n\theta\le2$ when $x\le\frac12$, and $1-\cos n\theta\le\frac12n^2\theta^2$ together with $\theta^2(\delta^2+\theta^2)^{-2}\le(\delta^2+\theta^2)^{-1}$ when $x\ge\frac12$,
		\begin{equation}\label{eq:d2Pint}
			\int_{-\pi}^{\pi}(1-\cos n\theta)\,|\partial_\theta^2P(x,\theta)|\,d\theta\le\begin{cases}640\pi x,&x\le\frac12,\\[2pt]\frac{5\pi^4n^2}2\displaystyle\int_\RR\frac{\delta\,d\theta}{\delta^2+\theta^2}=\frac{5\pi^5}2n^2,&x\ge\frac12.\end{cases}
		\end{equation}
		
		Now let $\omega$ be any weight and $j_\omega$ the function \eqref{eq:jomega} built from $\tilde\omega$. With $x=r^2$, \eqref{eq:Ctheta} bounds the integrand of \eqref{eq:jomega} by $C(\theta)\,\tilde\omega(r)\,r$, which is integrable since $\int_0^1\tilde\omega\,r\,dr=\frac12\|\omega\|_1$. So $j_\omega$ is well defined, continuous on $\TT\setminus\{0\}$ by dominated convergence, and even by \eqref{eq:d2Pseries}. Moreover \eqref{eq:d2Pint}, Tonelli's theorem and the inequality $1/r\le2r$ for $r^2\ge\frac12$ give
		$$\int_0^1\frac{\tilde\omega(r)}{r}\int_{-\pi}^{\pi}(1-\cos n\theta)\,|\partial_\theta^2P(r^2,\theta)|\,d\theta\,dr\le\big(640\pi+5\pi^5\big)n^2\int_0^1\tilde\omega(r)\,r\,dr<\infty .$$
		In particular $\int_{-\pi}^{\pi}(1-\cos n\theta)\,|j_\omega(\theta)|\,d\theta<\infty$ for every $n$, and Fubini's theorem applies. Two integrations by parts on the circle give $\int_{-\pi}^{\pi}(1-\cos n\theta)\,\partial_\theta^2P(x,\theta)\,d\theta=n^2\int_{-\pi}^{\pi}P(x,\theta)\cos n\theta\,d\theta=2\pi n^2x^{|n|}$, and with \eqref{eq:psi} we obtain the identity that links $j_\omega$ to the energy:
		\begin{equation}\label{eq:fourier}
			\psi_\omega(n)=\frac2\pi\int_{-\pi}^{\pi}(1-\cos n\theta)\,j_\omega(\theta)\,d\theta\qquad(n\in\ZZ),
		\end{equation}
		the integral being absolutely convergent. We also record how $j_\omega$ is recovered from $\psi_\omega$. For $\rho<1$, \eqref{eq:psi} and \eqref{eq:d2Pseries} with $x=\rho r^2$ give
		\begin{equation}\label{eq:abelsum}
			\sum_{n\in\ZZ}\psi_\omega(n)\rho^{|n|}e^{in\theta}=\int_0^1\frac{\tilde\omega(r)}{r}\,4\sum_{n\in\ZZ}n^2(\rho r^2)^{|n|}e^{in\theta}\,dr=-4\int_0^1\frac{\tilde\omega(r)}{r}\,\partial_\theta^2P(\rho r^2,\theta)\,dr,
		\end{equation}
		the interchange being legitimate because $\sum_nn^2\rho^{|n|}<\infty$ and $\int_0^1\tilde\omega\,r^{2|n|-1}dr\le\frac12\|\omega\|_1$ for $n\ne0$. As $\rho\uparrow1$ the last integrand converges pointwise and is dominated by $C(\theta)\,\tilde\omega(r)\,r$, so that, for $\theta\ne0$,
		\begin{equation}\label{eq:abel}
			j_\omega(\theta)=-\frac14\lim_{\rho\uparrow1}\sum_{n\in\ZZ}\psi_\omega(n)\rho^{|n|}e^{in\theta}.
		\end{equation}
		
		We can now prove that (ii) implies (i). It is convenient to argue with an arbitrary symmetric Borel measure $\tilde\nu\ge0$ on $\TT\setminus\{0\}$ such that
		\begin{equation}\label{eq:psinu}
			\psi_\omega(n)=\int|1-e^{in\theta}|^2\,\tilde\nu(d\theta)\qquad(n\in\ZZ).
		\end{equation}
		If $j_\omega\ge0$ on $(0,\pi]$, then $j_\omega\ge0$ on $\TT\setminus\{0\}$ because $j_\omega$ is even, and \eqref{eq:fourier}, with $|1-e^{in\theta}|^2=2(1-\cos n\theta)$, shows that $\tilde\nu(d\theta)=\frac1\pi j_\omega(\theta)\,d\theta$ satisfies \eqref{eq:psinu}. Let $\J$ be the image of $m\otimes\tilde\nu$ under $(\zeta,\theta)\mapsto(\zeta,e^{i\theta}\zeta)$. The exchange of $\zeta$ and $\lambda=e^{i\theta}\zeta$ corresponds to $(\zeta,\theta)\mapsto(e^{i\theta}\zeta,-\theta)$, which preserves $m\otimes\tilde\nu$, so $\J$ is symmetric; it is locally finite on $\TT^2\setminus\diag$ because $\int(1-\cos\theta)\,d\tilde\nu=\frac12\psi_\omega(1)<\infty$. Let $f$ be a Borel function with $\int|f|^2dm<\infty$. The function $\zeta\mapsto f(e^{i\theta}\zeta)$ has Fourier coefficients $\hat f(n)e^{in\theta}$, so Tonelli's theorem, Parseval's identity, \eqref{eq:psinu} and \eqref{eq:parseval} give
		\begin{equation}\label{eq:QJradial}
			Q_\J(f)=\int\tilde\nu(d\theta)\int_\TT|f(\zeta)-f(e^{i\theta}\zeta)|^2dm(\zeta)=\int\tilde\nu(d\theta)\sum_n|\hat f(n)|^2|1-e^{in\theta}|^2=\sum_n\psi_\omega(n)|\hat f(n)|^2=Q_\omega(f).
		\end{equation}
		Hence $\omega\in\WLtwo$. With $\tilde\nu=\frac1\pi j_\omega\,d\theta$ and $\theta=\arg\bar\zeta\lambda$ we have $dm(\lambda)=d\theta/2\pi$, so $\frac1\pi\,d\theta=2\,dm(\lambda)$ and $d\J=2j_\omega(\arg\bar\zeta\lambda)\,dm(\zeta)\,dm(\lambda)$. That (i) implies (i') is clear, since Lipschitz functions are Borel functions in $L^2$.
		
		We next prove that (i') implies (ii). The argument does not use radiality, and gives more. Let $\omega$ be any weight and let $\J$ be a symmetric positive Borel measure on $\TT^2\setminus\diag$ such that $Q_\omega(e_n)=Q_\J(e_n)$ for every $n\ge1$; this holds when $\omega$ is a Douglas weight and $\J=\J_\omega$. Put $\theta(\zeta,\lambda)=\arg\bar\zeta\lambda$, a continuous map from $\TT^2\setminus\diag$ to $\TT\setminus\{0\}$, and let $\theta_\#\J$ be the image of $\J$ under this map, that is, the Borel measure on $\TT\setminus\{0\}$ defined by $\theta_\#\J(E)=\J\big(\{(\zeta,\lambda):\arg\bar\zeta\lambda\in E\}\big)$, so that $\int\chi\,d(\theta_\#\J)=\iint\chi(\arg\bar\zeta\lambda)\,d\J(\zeta,\lambda)$ for every Borel function $\chi\ge0$ on $\TT\setminus\{0\}$. We claim that
		\begin{equation}\label{eq:necessary}
			\theta_\#\J=\frac1\pi\,j_\omega(\theta)\,d\theta\qquad\text{on }\TT\setminus\{0\},
		\end{equation}
		where $j_\omega$ is built from $\tilde\omega$; in particular $j_\omega\ge0$, which is (ii). Let $\nu_1=\theta_\#\J$. Since $\J$ is symmetric and the exchange of $\zeta$ and $\lambda$ changes $\theta$ into $-\theta$, $\nu_1$ is symmetric. As $|e_n(\zeta)-e_n(\lambda)|^2=|1-e^{in\theta}|^2=2(1-\cos n\theta)$, \eqref{eq:psi} gives $\psi_\omega(n)=Q_\J(e_n)=2\int(1-\cos n\theta)\,d\nu_1$ for $n\ge1$. The measure $\nu_1$ is locally finite on $\TT\setminus\{0\}$, because $\J(\{|\zeta-\lambda|>\varepsilon\})\le\varepsilon^{-2}\iint|\zeta-\lambda|^2d\J=\varepsilon^{-2}\psi_\omega(1)$. Put $\mu_1=2\nu_1$ and $\mu_2=\frac2\pi j_\omega\,d\theta$. By what precedes and by \eqref{eq:fourier}, $\int(1-\cos n\theta)\,d\mu_1=\int(1-\cos n\theta)\,d\mu_2=\psi_\omega(n)$ for every $n\ge0$, and all these integrals converge absolutely. Taking second differences in $n$ and using $2\cos n\theta-\cos(n+1)\theta-\cos(n-1)\theta=2\cos n\theta\,(1-\cos\theta)$, we get $\int\cos n\theta\,(1-\cos\theta)\,d\mu_1=\int\cos n\theta\,(1-\cos\theta)\,d\mu_2$ for $n\ge1$, and the case $n=0$ is the identity for $n=1$. The measures $(1-\cos\theta)\mu_1$ and $(1-\cos\theta)\mu_2$ are finite, the second one possibly signed, and both are symmetric; they have the same Fourier coefficients, so they coincide. As $1-\cos\theta>0$ on $\TT\setminus\{0\}$, $\mu_1=\mu_2$ there, which is \eqref{eq:necessary}, and the continuity of $j_\omega$ gives $j_\omega\ge0$ on $(0,\pi]$.
		
		It remains to identify $\J_\omega$, and for this we show that a weight has at most one symmetric Douglas measure. Let $\omega$ be any weight and let $\J$, $\J'$ be symmetric positive Borel measures on $\TT^2\setminus\diag$ with $Q_\omega(f)=Q_\J(f)=Q_{\J'}(f)$ for every Lipschitz function $f$. These forms are finite on Lipschitz functions, since $Q_\J(f)\le\mathrm{Lip}(f)^2\iint|\zeta-\lambda|^2d\J=\mathrm{Lip}(f)^2\,Q_\omega(e_1)$. Let $g,h\ge0$ be real Lipschitz functions with disjoint supports, and put $a=g(\zeta)-g(\lambda)$, $b=h(\zeta)-h(\lambda)$. Since $gh\equiv0$, $ab=-g(\zeta)h(\lambda)-g(\lambda)h(\zeta)$, and polarization gives $Q_\omega(g,h)=\frac14\big(Q_\omega(g+h)-Q_\omega(g-h)\big)=\frac14\iint\big(|a+b|^2-|a-b|^2\big)\,d\J=\iint ab\,d\J$, that is, by symmetry of $\J$,
		\begin{equation}\label{eq:offdiag}
			Q_\omega(g,h)=-\iint\big(g(\zeta)h(\lambda)+g(\lambda)h(\zeta)\big)\,d\J=-2\iint g(\zeta)h(\lambda)\,d\J(\zeta,\lambda),
		\end{equation}
		and the same holds for $\J'$. Let $I$, $I'$ be open arcs of length less than $\frac\pi2$ with disjoint closures, at distance $d>0$, and let $g_k=\min(1,k\,\mathrm{dist}(\cdot,\TT\setminus I))$ and $h_k=\min(1,k\,\mathrm{dist}(\cdot,\TT\setminus I'))$. These are nonnegative Lipschitz functions supported in $I$ and $I'$, and $g_k\uparrow\mathbf 1_I$, $h_k\uparrow\mathbf 1_{I'}$. By \eqref{eq:offdiag} for $\J$ and $\J'$ and monotone convergence, $\J(I\times I')=\J'(I\times I')$, and this number is finite, being at most $\J(\{|\zeta-\lambda|\ge d\})\le d^{-2}Q_\omega(e_1)$. Two arcs of length less than $\frac\pi2$ meet in an arc or not at all, so these rectangles form a family closed under finite intersections; it generates the Borel sets of $\TT^2\setminus\diag$, and countably many of its members, with rational endpoints, cover $\TT^2\setminus\diag$. Hence $\J=\J'$. For a radial Douglas weight, $\J_\omega$ is therefore the measure $2j_\omega(\arg\bar\zeta\lambda)\,dm\,dm$ constructed above. This completes the proof.
	\end{proof}
	
	For radial weights, the two sides of the Douglas formula have the same domain with the notation of Section 1 
	
	\begin{corollary}\label{cor:domains}
		Let $\omega$ be a radial Douglas weight. Then $D^h(\omega)=F_\omega=D(\J_\omega)$, and this space consists of the functions $f\in L^2(\TT)$ with $\sum_{n\in\ZZ}\psi_\omega(n)|\hat f(n)|^2<\infty$.
	\end{corollary}
	\begin{proof}
		By \eqref{eq:parseval}, $D^h(\omega)$ is the space described in the statement. If $f$ belongs to it, the partial Fourier sums $S_Nf$ are smooth, converge to $f$ in $L^2$, and $Q_\omega(f-S_Nf)=\sum_{|n|>N}\psi_\omega(n)|\hat f(n)|^2\to0$; so $f\in F_\omega$. Conversely, if smooth functions $f_k$ converge to $f$ in $L^2$ and form a Cauchy sequence for $Q_\omega$, then $P[f_k]\to P[f]$ locally uniformly in $\DD$ with all derivatives, and Fatou's lemma gives $Q_\omega(f)\le\liminf_kQ_\omega(f_k)<\infty$; so $F_\omega\subset D^h(\omega)$. By Theorem~\ref{thm:A}, $\J_\omega$ is the image of $m\otimes\tilde\nu$, with $\tilde\nu=\frac1\pi j_\omega\,d\theta$, under $(\zeta,\theta)\mapsto(\zeta,e^{i\theta}\zeta)$. For an $m$-null set $N$, Tonelli's theorem gives $\J_\omega(N\times\TT)=\int m(N)\,\tilde\nu(d\theta)=0$, and similarly $\J_\omega(\TT\times N)=0$ by symmetry, so $Q_{\J_\omega}(f)$ depends only on the class of $f$ in $L^2$. By \eqref{eq:QJradial}, $Q_{\J_\omega}(f)=Q_\omega(f)$ for every such class, hence $D(\J_\omega)=D^h(\omega)$.
	\end{proof}
	
	\section{Proof of Theorem~\ref{thm:B}}\label{sec:main}
	
	\paragraph{\textbf{Outline of the proof}}
	We consider the distribution $K_\omega$ defined below.
	Suppose first that $\omega$ is a Douglas weight with associated
	measure $\J$. By polarizing \eqref{eq:intro-def} (see
	\eqref{eq:offdiag}), we obtain
	\[
	Q_\omega(g,h)
	=-\iint_{\TT^2}
	\bigl(g(\zeta)h(\lambda)+g(\lambda)h(\zeta)\bigr)
	\,d\J(\zeta,\lambda)
	\]
	for all nonnegative functions $g,h$ with disjoint supports.
	It follows that $K_\omega$ is nonpositive off the diagonal
	and coincides there with $-2\J$. Moreover, Lemma~\ref{lem:elem}
	yields
	\[
	\langle K_\omega,\varphi\rangle
	=\frac12\int_\DD \omega\,\Delta(\mathcal P\varphi)\,dA,
	\]
	so that the corresponding sign condition is precisely
	\eqref{eq:weaksuper} in Definition~\ref{def:weaksuper}.
	
	Conversely, assume that $\omega$ satisfies \eqref{eq:weaksuper}.
	Then the restriction
	\[
	\J=-\frac12 K_\omega\big|_{\TT^2\setminus\diag}
	\]
	defines a positive measure on $\TT^2\setminus\diag$.
	We compare the quadratic forms $Q_\omega$ and $Q_\J$.
	Since they agree on pairs of functions with disjoint supports,
	their difference is a distribution supported on the diagonal.
	Testing $\J$ against functions depending on
	$\arg(\bar\zeta\lambda)$, we obtain
	\[
	\iint_{\TT^2}|\zeta-\lambda|^2\,d\J(\zeta,\lambda)<\infty.
	\]
	This integrability condition, together with the growth estimates
	for the matrix $M_\omega$ and the structure theorem for
	distributions supported on the diagonal, implies that the
	difference between $Q_\omega$ and $Q_\J$ vanishes.
	Finally, Lemma~\ref{lem:approx} extends the resulting identity
	to Lipschitz functions.
	
	Theorem~\ref{thm:B} proves the equivalence of (i) and (iii) below
	and identifies the corresponding measure. Condition~(ii) expresses
	\eqref{eq:weaksuper} in terms of the bilinear form $Q_\omega$
	and is the form used in the proof.
	
	\begin{theorem}\label{thm:main}
		Let $\omega\ge0$ be measurable with $\int_\DD\omega\,dA<\infty$. The following are equivalent.
		\begin{enumerate}
			\item[(i)] $\omega\in\WLip$, that is, $\omega$ is a Douglas weight.
			\item[(ii)] $Q_\omega(g,h)\le0$ for all nonnegative $g,h\in C^\infty(\TT)$ with disjoint supports.
			\item[(iii)] $\omega$ satisfies condition~\eqref{eq:weaksuper}.
		\end{enumerate}
		In this case the symmetric measure in (i) is unique and is given by
		$$\iint\varphi\,d\J_\omega=-\frac14\int_\DD\omega\,\Delta\big(\mathcal P\varphi\big)\,dA\qquad\big(\varphi\in C_c^\infty(\TT^2\setminus\diag)\big),$$
		and $\iint|\zeta-\lambda|^2d\J_\omega=\psi_\omega(1)=2\|\omega\|_1$.
	\end{theorem}
	
	We begin with four lemmas. The first reduces the Douglas formula
	for Lipschitz functions to the case of trigonometric polynomials.
	
	\begin{lemma}\label{lem:approx}
		Let $\J$ be a symmetric positive Borel measure on $\TT^2\setminus\diag$ with $\iint|\zeta-\lambda|^2d\J<\infty$. If $Q_\omega(p)=Q_\J(p)$ for every trigonometric polynomial $p$, then $Q_\omega(f)=Q_\J(f)<\infty$ for every $f\in\mathrm{Lip}(\TT)$.
	\end{lemma}
	\begin{proof}
		Let $f\in\mathrm{Lip}(\TT)$ and let $f_k$ be its Fej\'er means. They are trigonometric polynomials, $f_k\to f$ uniformly, and $\mathrm{Lip}(f_k)\le\mathrm{Lip}(f)$, because $f_k$ is an average of rotations of $f$. Hence $|(f_k-f)(\zeta)-(f_k-f)(\lambda)|^2\le4\,\mathrm{Lip}(f)^2|\zeta-\lambda|^2$, which is $\J$-integrable, and dominated convergence gives $Q_\J(f_k-f)\to0$. In particular $Q_\omega(f_k-f_l)=Q_\J(f_k-f_l)\to0$ as $k,l\to\infty$. Since $P[f_l]\to P[f]$ locally uniformly in $\DD$ together with all derivatives, Fatou's lemma gives
		$$Q_\omega(f-f_k)\le\liminf_{l\to\infty}Q_\omega(f_l-f_k)=\liminf_{l\to\infty}Q_\J(f_l-f_k)=Q_\J(f-f_k),$$
		which tends to $0$ as $k\to\infty$. As $Q_\omega^{1/2}$ and $Q_\J^{1/2}$ are seminorms, $Q_\omega(f)=\lim_kQ_\omega(f_k)=\lim_kQ_\J(f_k)=Q_\J(f)$.
	\end{proof}
	
	For $s\ge2$ and $j\in\ZZ$ put $\Omega_j(s)=\int_0^1r^{s-1}\omega_j(r)\,dr$ . We have $|\Omega_j(s)|\le\int_0^1r\,\tilde\omega\,dr=\frac12\|\omega\|_1$, and $\Omega_j(s)\to0$ as $s\to\infty$ by dominated convergence.
	
	\begin{lemma}\label{lem:matrix}
		For $n,k\in\ZZ$,
		$$M_\omega(n,k)=\E_\omega(e_n,e_k)=4|n||k|\,\mathbf 1_{\{nk>0\}}\,\Omega_{k-n}(|n|+|k|).$$
		In particular $M_\omega(n,k)=0$ when $nk\le0$, $M_\omega(n,n)=\psi_\omega(n)$, $M_\omega(-n,-k)=\overline{M_\omega(n,k)}$ and $|M_\omega(n,k)|\le2|n||k|\,\|\omega\|_1$. For the bilinear form, $Q_\omega(e_n,e_m)=M_\omega(n,-m)$.
	\end{lemma}
	\begin{proof}
		Let $n,k\ge1$, $u=z^n$ and $v=z^k$. Then $\nabla u\cdot\overline{\nabla v}=u_x\bar v_x+u_y\bar v_y=2nk\,z^{n-1}\bar z^{k-1}=2nk\,r^{n+k-2}e^{i(n-k)t}$ for $z=re^{it}$, and \eqref{eq:polar} gives
		$$\E_\omega(e_n,e_k)=4nk\int_0^1r^{n+k-1}\Big(\frac1{2\pi}\int_0^{2\pi}e^{-i(k-n)t}\omega(re^{it})\,dt\Big)dr=4nk\,\Omega_{k-n}(n+k).$$
		If $n>0>k$, then $u=z^n$, $v=\bar z^{|k|}$, $\bar v$ is holomorphic and $u_x\bar v_x+u_y\bar v_y=n|k|\,z^{n-1}z^{|k|-1}(1+i\cdot i)=0$. The case $n,k<0$ follows by conjugation, and the bound from $|\Omega_j(s)|\le\frac12\|\omega\|_1$ for $s\ge2$.
	\end{proof}
	
	By Lemma~\ref{lem:matrix}, $(g,h)\mapsto Q_\omega(g,h)$ is a continuous bilinear form on $C^\infty(\TT)$, because the Fourier coefficients of smooth functions decay rapidly. It therefore defines a symmetric distribution $K_\omega\in\mathcal D'(\TT^2)$ by $\langle K_\omega,g\otimes h\rangle=Q_\omega(g,h)$, the pairing being taken with respect to $m\otimes m$. In Fourier coordinates, with $\zeta=e^{i\alpha}$ and $\lambda=e^{i\beta}$,
	$$K_\omega(\alpha,\beta)=\sum_{n,k\in\ZZ}M_\omega(n,k)\,e^{-in\alpha+ik\beta}.$$
	
	On functions of the angle $\arg\bar\zeta\lambda$, the distribution $K_\omega$ depends only on $j_\omega$; this is part (b) of the next lemma, which also expresses $K_\omega$ through the Poisson kernels.
	
	\begin{lemma}\label{lem:elem}
		\textup{(a)} For $z\in\DD$ and $\zeta,\lambda\in\TT$ put $\kappa_z(\zeta,\lambda)=\nabla P_\zeta(z)\cdot\nabla P_\lambda(z)$. Then
		\begin{equation}\label{eq:kappaz}
			\kappa_z(\zeta,\lambda)=4\,\Re\frac{\zeta\bar\lambda}{(\zeta-z)^2\,\overline{(\lambda-z)}^{\,2}},\qquad K_\omega=\int_\DD\kappa_z\,\omega(z)\,dA(z)\ \text{ in }\mathcal D'(\TT^2).
		\end{equation}
		For $\varphi\in C^\infty(\TT^2)$, the function $z\mapsto\langle\kappa_z,\varphi\rangle=\frac12\Delta\big(\mathcal P\varphi\big)(z)$ is bounded on $\DD$, where $\mathcal P\varphi$ is defined by \eqref{eq:intro-Pphi}, and $\langle K_\omega,\varphi\rangle=\frac12\int_\DD\omega\,\Delta(\mathcal P\varphi)\,dA$.
		
		\textup{(b)} Let $\chi\in C_c^\infty(\TT\setminus\{0\})$ and $\varphi_\chi(\zeta,\lambda)=\chi(\arg\bar\zeta\lambda)$, so that $\varphi_\chi\in C_c^\infty(\TT^2\setminus\diag)$. Then
		$$\langle K_\omega,\varphi_\chi\rangle=-4\int_{-\pi}^{\pi}\chi(\theta)\,j_\omega(\theta)\,\frac{d\theta}{2\pi}.$$
	\end{lemma}
	\begin{proof}
		\textup{(a)} We have $P_\zeta=\Re F_\zeta$ with $F_\zeta(z)=\frac{\zeta+z}{\zeta-z}$ and $F_\zeta'(z)=\frac{2\zeta}{(\zeta-z)^2}$. If $u=\Re F$ and $v=\Re G$ with $F$, $G$ holomorphic, then $u_x=\Re F'$, $u_y=-\Im F'$ and $\nabla u\cdot\nabla v=\Re(F'\overline{G'})$, which gives the formula for $\kappa_z$. Let $\varphi=\sum\hat\varphi(n,k)\,e_n\otimes e_k$ be smooth. Since $\langle\kappa_z,e_n\otimes e_k\rangle=\nabla P[e_n](z)\cdot\nabla P[e_k](z)$ and $|\nabla P[e_n](z)|=\sqrt2\,|n|\,|z|^{|n|-1}$, we get $|\langle\kappa_z,\varphi\rangle|\le2\sum|\hat\varphi(n,k)||n||k|$ for all $z\in\DD$. Fubini's theorem then gives $\int_\DD\omega\,\langle\kappa_z,\varphi\rangle\,dA=\sum\hat\varphi(n,k)\,Q_\omega(e_n,e_k)=\langle K_\omega,\varphi\rangle$. Finally, $\Delta_z(P_\zeta P_\lambda)=2\nabla P_\zeta\cdot\nabla P_\lambda$ because $P_\zeta$ and $P_\lambda$ are harmonic, so $\langle\kappa_z,\varphi\rangle=\frac12\Delta\mathcal P\varphi(z)$.
		
		\textup{(b)} Write $\chi(\theta)=\sum_n\hat\chi(n)e^{in\theta}$ with rapidly decreasing coefficients. Then $\varphi_\chi=\sum_n\hat\chi(n)\,e_{-n}\otimes e_n$ in $C^\infty(\TT^2)$, and $Q_\omega(e_{-n},e_n)=M_\omega(-n,-n)=\psi_\omega(n)$. As $\psi_\omega(n)\le2n^2\|\omega\|_1$,
		$$\langle K_\omega,\varphi_\chi\rangle=\sum_{n\in\ZZ}\hat\chi(n)\,\psi_\omega(n)=\lim_{\rho\uparrow1}\int_{-\pi}^{\pi}\chi(\theta)\sum_{n\in\ZZ}\psi_\omega(n)\rho^{|n|}e^{-in\theta}\,\frac{d\theta}{2\pi}.$$
		By \eqref{eq:abelsum} and \eqref{eq:Ctheta}, the inner sum is bounded by $4C(\theta)\int_0^1\tilde\omega(r)\,r\,dr$, with $C(\theta)$ bounded on the support of $\chi$, and it converges to $-4j_\omega(\theta)$. Dominated convergence gives the result.
	\end{proof}
	
	The last lemma reduces test functions on $\TT^2\setminus\diag$ to tensor products $g\otimes h$ with disjoint supports.
	
	\begin{lemma}\label{lem:tensor}
		Every nonnegative $\varphi\in C_c^\infty(\TT^2\setminus\diag)$ is a finite sum of nonnegative functions, each of which is a limit in $C^k$, for every $k$, of finite sums $\sum_ig_i\otimes h_i$ with $g_i,h_i\ge0$ smooth and $\supp g_i\cap\supp h_i=\emptyset$. The same holds without positivity for real $\varphi$.
	\end{lemma}
	\begin{proof}
		Since $\supp\varphi$ is at positive distance from $\diag$, it is covered by finitely many open rectangles $I\times I'$ with $\overline I\cap\overline{I'}=\emptyset$. A smooth partition of unity splits $\varphi$ into nonnegative pieces $\psi$, each supported in one rectangle $I\times I'$. Convolving $\psi$ with $\eta_\varepsilon\otimes\eta_\varepsilon$, where $\eta_\varepsilon\ge0$ is a small bump, keeps the support in $I\times I'$, and replacing the convolution integral by Riemann sums gives nonnegative combinations of products of bumps supported in $I$ and $I'$, which converge to $\psi$ in $C^k$.
	\end{proof}
	
	\begin{proof}[Proof of Theorem~\ref{thm:main}]
		By Lemma~\ref{lem:matrix}, $K_\omega$ is a distribution on the compact manifold $\TT^2$, so it has finite order: there are $k\ge0$ and $C>0$ with $|\langle K_\omega,\varphi\rangle|\le C\|\varphi\|_{C^k}$ for every $\varphi\in C^\infty(\TT^2)$. By Lemma~\ref{lem:elem}(a), $\langle K_\omega,\varphi\rangle=\frac12\int_\DD\omega\,\Delta(\mathcal P\varphi)\,dA$, so condition (iii) says exactly that $\langle K_\omega,\varphi\rangle\le0$ for every nonnegative $\varphi\in C_c^\infty(\TT^2\setminus\diag)$, that is, that $K_\omega$ is a nonpositive distribution on the open set $\TT^2\setminus\diag$. If $g,h\ge0$ are smooth with disjoint supports, then $g\otimes h$ is such a function and $\langle K_\omega,g\otimes h\rangle=Q_\omega(g,h)$, so (iii) implies (ii). Conversely, (ii) says that $\langle K_\omega,\cdot\rangle\le0$ on finite sums $\sum_ig_i\otimes h_i$ with $g_i,h_i\ge0$ smooth and $\supp g_i\cap\supp h_i=\emptyset$; by the continuity of $K_\omega$ for the $C^k$ norm, the same holds for their $C^k$ limits, and Lemma~\ref{lem:tensor} gives $\langle K_\omega,\varphi\rangle\le0$ for every nonnegative $\varphi\in C_c^\infty(\TT^2\setminus\diag)$. Thus (ii) and (iii) are equivalent. If (i) holds with a measure $\J$ and $g,h\ge0$ are smooth with disjoint supports, then $g$ and $h$ are Lipschitz, and \eqref{eq:offdiag} gives $Q_\omega(g,h)=-2\iint g(\zeta)h(\lambda)\,d\J\le0$; so (i) implies (ii).
		
		It remains to prove that (ii) implies (i), and we assume (ii) from now on. A positive distribution on an open set is a positive Radon measure \cite[Section~2.1]{hormander1983analysis}, so
		\begin{equation}\label{eq:Jdef}
			\J=-\tfrac12K_\omega|_{\TT^2\setminus\diag}
		\end{equation}
		is a positive Radon measure on $\TT^2\setminus\diag$. It is symmetric because $Q_\omega(g,h)=Q_\omega(h,g)$, and by construction
		\begin{equation}\label{eq:QomegaJ}
			Q_\omega(g,h)=\langle K_\omega,g\otimes h\rangle=-2\iint g(\zeta)h(\lambda)\,d\J(\zeta,\lambda)\qquad(g,h\in C^\infty(\TT),\ \supp g\cap\supp h=\emptyset).
		\end{equation}
		We show that $\J$ is the Douglas measure of $\omega$.
		
		The first point is that $\J$ has a finite second moment, and this is where $j_\omega$ enters. Let $\chi\ge0$ belong to $C_c^\infty(\TT\setminus\{0\})$. By \eqref{eq:Jdef} and Lemma~\ref{lem:elem}(b),
		$$\iint\chi(\arg\bar\zeta\lambda)\,d\J(\zeta,\lambda)=-\frac12\langle K_\omega,\varphi_\chi\rangle=\frac1\pi\int_{-\pi}^{\pi}\chi(\theta)\,j_\omega(\theta)\,d\theta .$$
		The left-hand side is nonnegative, and $j_\omega$ is continuous on $\TT\setminus\{0\}$, so $j_\omega\ge0$ there. The image of $\J$ under $(\zeta,\lambda)\mapsto\arg\bar\zeta\lambda$ and the measure $\frac1\pi j_\omega\,d\theta$ are therefore two positive Radon measures on $\TT\setminus\{0\}$ that agree on $C_c^\infty(\TT\setminus\{0\})$; they are equal, and the identity above holds for every Borel function $\chi\ge0$. Since $|\zeta-\lambda|^2=|1-e^{i\theta}|^2=2(1-\cos\theta)$ for $\theta=\arg\bar\zeta\lambda$, the choice $\chi(\theta)=2(1-\cos\theta)$ and \eqref{eq:fourier} with $n=1$ give
		\begin{equation}\label{eq:secondmoment}
			\iint|\zeta-\lambda|^2\,d\J=\frac2\pi\int_{-\pi}^{\pi}(1-\cos\theta)\,j_\omega(\theta)\,d\theta=\psi_\omega(1)=2\|\omega\|_1 .
		\end{equation}
		Consequently $Q_\J(g,h)$ is defined for Lipschitz functions $g$, $h$, and $|Q_\J(g,h)|\le\mathrm{Lip}(g)\,\mathrm{Lip}(h)\,\psi_\omega(1)$.
		
		We now compare $Q_\omega$ and $Q_\J$ on smooth functions. By Lemma~\ref{lem:matrix} and the bound just obtained, $R=Q_\omega-Q_\J$ is a continuous bilinear form on $C^\infty(\TT)\times C^\infty(\TT)$. It vanishes on pairs of functions with disjoint supports. Indeed, for such $g$ and $h$ the computation that led to \eqref{eq:offdiag}, which uses only the symmetry of $\J$ and $gh\equiv0$, gives $Q_\J(g,h)=-2\iint g(\zeta)h(\lambda)\,d\J$, and this is $Q_\omega(g,h)$ by \eqref{eq:QomegaJ}. Moreover $R(1,\cdot)=R(\cdot,1)=0$, because both forms vanish as soon as one argument is constant. By the Schwartz kernel theorem \cite[Theorem~5.2.1]{hormander1983analysis} there is $\rho\in\mathcal D'(\TT^2)$ with $\langle\rho,g\otimes h\rangle=R(g,h)$; it has some finite order $N$, and since it vanishes on the sums $\sum_ig_i\otimes h_i$ of Lemma~\ref{lem:tensor}, it vanishes on $C_c^\infty(\TT^2\setminus\diag)$. Thus $\rho$ is supported on the diagonal.
		
		A distribution supported on the diagonal is killed by every function that vanishes there to a sufficiently high order, and this forces a rigid structure on the Fourier matrix $B(n,m)=\langle\rho,e_n\otimes e_m\rangle=R(e_n,e_m)$. The function $(\zeta-\lambda)^{N+1}e_n(\zeta)e_m(\lambda)$ vanishes to order $N+1$ on $\diag$, so $\langle\rho,(\zeta-\lambda)^{N+1}e_n\otimes e_m\rangle=0$ \cite[Theorem~2.3.3]{hormander1983analysis}, and the binomial formula gives
		$$\sum_{a=0}^{N+1}(-1)^{N+1-a}\binom{N+1}a\,B(n+a,\,m+N+1-a)=0\qquad(n,m\in\ZZ).$$
		In each term the two indices add up to $n+m+N+1$. Hence, for fixed $s\in\ZZ$, the function $b_s(n)=B(n,s-n)$ has vanishing $(N+1)$-st difference, and is a polynomial in $n$ of degree at most $N$. Its growth is controlled by the two forms: Lemma~\ref{lem:matrix} gives $|Q_\omega(e_n,e_m)|=|M_\omega(n,-m)|\le2|n||m|\,\|\omega\|_1$, and $|\zeta^n-\lambda^n|\le|n||\zeta-\lambda|$ together with \eqref{eq:secondmoment} gives $|Q_\J(e_n,e_m)|\le|n||m|\,\psi_\omega(1)$. So $b_s(n)=O(n^2)$ and $b_s$ has degree at most $2$. Moreover $b_s(0)=b_s(s)=0$, since $R(1,\cdot)=R(\cdot,1)=0$. If $s\ne0$, a polynomial of degree at most $2$ vanishing at $0$ and $s$ is $c_s\,n(s-n)$; if $s=0$, $b_0$ is even by the symmetry of $B$ and vanishes at $0$, so $b_0(n)=c_0n^2=-c_0\,n(0-n)$. In both cases there is $\gamma_s\in\CC$ with
		$$B(n,m)=\gamma_{n+m}\,nm\qquad(n,m\in\ZZ).$$
		
		The coefficients $\gamma_s$ vanish because both forms are $o(nm)$ along each anti-diagonal. Fix $s$, let $n\to+\infty$ and put $m=s-n$, so that $-m=n-s\ge1$ for $n$ large. By Lemma~\ref{lem:matrix}, $Q_\omega(e_n,e_m)=M_\omega(n,n-s)=4n(n-s)\,\Omega_{-s}(2n-s)$, hence $Q_\omega(e_n,e_m)/(nm)=-4\,\Omega_{-s}(2n-s)\to0$. On the other hand,
		$$\frac{Q_\J(e_n,e_m)}{nm}=\iint\frac{\zeta^n-\lambda^n}n\cdot\frac{\zeta^m-\lambda^m}m\,d\J\longrightarrow0$$
		by dominated convergence: the integrand is bounded by $|\zeta-\lambda|^2$, which is $\J$-integrable by \eqref{eq:secondmoment}, and it is at most $4/(n|m|)$ in absolute value. Therefore $\gamma_s=\lim B(n,s-n)/(n(s-n))=0$ for every $s$, so $B\equiv0$. A distribution on $\TT^2$ with vanishing Fourier coefficients is zero, hence $\rho=0$ and $Q_\omega(g,h)=Q_\J(g,h)$ for all $g,h\in C^\infty(\TT)$.
		
		For a trigonometric polynomial $f$ this gives $Q_\omega(f)=Q_\omega(f,\bar f)=Q_\J(f,\bar f)=Q_\J(f)$. By \eqref{eq:secondmoment}, Lemma~\ref{lem:approx} applies, and $Q_\omega(f)=Q_\J(f)$ for every $f\in\mathrm{Lip}(\TT)$. This is (i), with the measure \eqref{eq:Jdef}.
		
		Finally, the symmetric measure in (i) is unique, by the argument at the end of Section~\ref{sec:radial}, so $\J_\omega$ is the measure \eqref{eq:Jdef}. For $\varphi\in C_c^\infty(\TT^2\setminus\diag)$, Lemma~\ref{lem:elem}(a) gives $\iint\varphi\,d\J_\omega=-\frac12\langle K_\omega,\varphi\rangle=-\frac14\int_\DD\omega\,\Delta(\mathcal P\varphi)\,dA$, and $\iint|\zeta-\lambda|^2d\J_\omega=2\|\omega\|_1$ is \eqref{eq:secondmoment}.
	\end{proof}
	
	\begin{remark}\label{rem:weak}
		For a nonnegative $\varphi\in C_c^\infty(\TT^2\setminus\diag)$, the function $\mathcal P\varphi$ is nonnegative and smooth in $\DD$ and tends to $0$ at every point of $\TT$, but its support is not compact in $\DD$. An integration by parts in $\int_\DD\omega\,\Delta(\mathcal P\varphi)\,dA$ therefore produces a boundary term, and \eqref{eq:weaksuper} is not a restriction of the inequality $\Delta\omega\le0$ to fewer test functions. For $\omega\equiv1$ we have $\Delta\omega=0$, whereas $\int_\DD\Delta(\mathcal P\varphi)\,dA=2\langle K_1,\varphi\rangle=-8\iint\varphi(\zeta,\lambda)\,|\zeta-\lambda|^{-2}\,dm\,dm$ is negative when $\varphi\ne0$ (Example~\ref{ex:j1}). For a bounded $C^2$ weight with bounded Laplacian, \eqref{eq:weaksuper} is the pointwise inequality \eqref{eq:pointwise} of Proposition~\ref{prop:smooth}, in which a positive Laplacian is allowed if the boundary values of $\omega$ compensate it; the weight $|z|^2$, with $\Delta|z|^2=4$, is a Douglas weight (Example~\ref{prop:subharm}).
	\end{remark}
	
	\begin{remark}\label{rem:holomorphic}
		If \eqref{eq:intro-def} is only required for boundary values of holomorphic functions, the measure is no longer unique, even for the constant weight. For $c\in\RR$ with $|c|\le\frac12$ put
		$$\J^{c}=\Big(\frac2{|\zeta-\lambda|^2}+c\,\Re(\zeta\lambda)\Big)\,dm(\zeta)\,dm(\lambda),$$
		a positive symmetric measure since $2|\zeta-\lambda|^{-2}\ge\frac12$. For a holomorphic polynomial $f$, the integrals $\int_\TT f(\zeta)\,\zeta\,dm(\zeta)$ and $\int_\TT\zeta\,dm(\zeta)$ vanish, and each term in the expansion of $|f(\zeta)-f(\lambda)|^2(\zeta\lambda+\bar\zeta\bar\lambda)$ has one of these integrals or its conjugate, in the variable $\zeta$ or $\lambda$, as a factor. Hence $\iint|f(\zeta)-f(\lambda)|^2\,\Re(\zeta\lambda)\,dm\,dm=0$ and $Q_{\J^c}(f)=Q_{\J^0}(f)=\E_1(f)$, although $\J^c\ne\J^0=\J_1$ when $c\ne0$. For radial weights the two requirements define the same class (Remark~\ref{rem:holoradial}); we do not know whether this holds in general.
	\end{remark}
	
	\section{Further results}\label{sec:further}
	
	We derive some consequences of Theorems~\ref{thm:A} and
	\ref{thm:B} that hold independently of the choice of weight.
	These consequences will be used in the study of examples in
	Section~\ref{sec:examples}.
	
	We begin by computing $j_\omega$ for two specific weights.
	The resulting formulas will be used throughout this section
	and in Section~\ref{sec:examples}.
	
	\begin{example}\label{ex:j1}
		For $\omega\equiv1$ we have $\psi_1(n)=2|n|$ and
		$$j_1(\theta)=\frac1{4\sin^2(\theta/2)}\qquad(\theta\ne0).$$
		With $\theta=\arg\bar\zeta\lambda$ we have $|\zeta-\lambda|^2=4\sin^2(\theta/2)$, so $j_1(\arg\bar\zeta\lambda)=|\zeta-\lambda|^{-2}$ and $\J_1=2|\zeta-\lambda|^{-2}dm\,dm$, in accordance with Proposition~\ref{prop:superharmonic} for $\mu=0$, $\nu=m$. The factor $2$ with respect to \eqref{eq:intro-douglas} comes from $\E_1(f)=2D(F)$.
	\end{example}
	
	For the weights considered in this paper, \eqref{eq:psi} gives, after the substitution $t=r^2$ in the last case,
	\begin{equation}\label{eq:psiexamples}
		\psi_1(n)=2|n|,\qquad\psi_{|z|^{2\gamma}}(n)=\frac{2n^2}{|n|+\gamma}\quad(\gamma\ge0),\qquad\psi_{(1-|z|^2)^s}(n)=2n^2B(|n|,s+1)\quad(s>-1),
	\end{equation}
	where $B$ is the Beta function. Thus $\psi_{1-|z|^2}(n)=2|n|/(|n|+1)$, and $2n^2B(|n|,s+1)\sim2\Gamma(s+1)|n|^{1-s}$ as $|n|\to\infty$.
	
	We next compute the kernel associated with the weight $1-|z|^2$.
	
	\begin{example}\label{lem:AdA}
		For $0<\theta<2\pi$, let $\mathsf A(\theta)=-1-2\cos\theta\,\log\big(2\sin\frac\theta2\big)+(\pi-\theta)\sin\theta$. We have
		\begin{enumerate}[label=\textup{(\alph*)}]
			\item $j_{1-|z|^2}(\theta)=\frac12\mathsf A(\theta)$;
			\item $\int_\DD P_\zeta P_\lambda\,dA=\mathsf A(\arg\bar\zeta\lambda)$ for $\zeta\ne\lambda$ in $\TT$; in particular $\mathsf A>0$;
			\item $\mathsf A(\theta)=2\log(1/\theta)-1+o(1)$ as $\theta\to0^+$;
			\item $4\sin^2(\theta/2)\,\mathsf A(\theta)\le2$.
		\end{enumerate}
	\end{example}
	
	Since $1-|z|^2=G[dA]$, Proposition~\ref{prop:superharmonic} gives $\J_{1-|z|^2}=A_{dA}\,dm\,dm$, and Theorem~\ref{thm:A} gives $\J_{1-|z|^2}=2j_{1-|z|^2}(\arg\bar\zeta\lambda)\,dm\,dm$; by (a) and (b) the two densities coincide. Numerically, $\sin^2(\theta/2)\mathsf A(\theta)$ is increasing on $(0,\pi]$, so its supremum is its value $2\log2-1\approx0.386$ at $\theta=\pi$; part (d) proves the weaker bound $\frac12$.
	
	\subsection{Necessary conditions}\label{subsec:necessary}
	
	Identity~\eqref{eq:necessary} does not require $\omega$ to be radial.
	Thus, if $\omega$ is a Douglas weight, the function $j_\omega$
	associated with its radialization $\tilde\omega$ is nonnegative.
	By Theorem~\ref{thm:A}, this also implies that $\tilde\omega$
	is a Douglas weight. The converse does not hold, as shown in
	Section~\ref{sec:examples}.
	
	Furthermore, the proof of \eqref{eq:necessary} uses only the
	functions $e_n$, $n\geq1$, which are boundary values of
	holomorphic polynomials. This observation leads to the following
	remark.
	
	\begin{remark}\label{rem:holoradial}
		Let $\omega$ be radial, and suppose that there is a symmetric positive measure $\J$ on $\TT^2\setminus\diag$ such that $\int_\DD|F'|^2\omega\,dA=\frac12Q_\J(f)$ for every holomorphic polynomial $F$ with boundary values $f$. Then $Q_\omega(e_n)=Q_\J(e_n)$ for $n\ge1$, so $j_\omega\ge0$ by \eqref{eq:necessary}, and $\omega$ is a Douglas weight by Theorem~\ref{thm:A}.
	\end{remark}
	
	The following consequence of \eqref{eq:necessary} holds for every weight.
	
	\begin{corollary}\label{cor:boundary}
		Let $\omega$ be a nonzero Douglas weight. Then $\psi_\omega(n)$ does not tend to $0$ as $n\to\infty$. In particular $\int_{\{\rho<|z|<1\}}\omega\,dA>0$ for every $\rho<1$, and no nonzero Douglas weight has compact support in $\DD$.
	\end{corollary}
	\begin{proof}
		Let $\tilde\nu=\theta_\#\J_\omega$ be the image of $\J_\omega$ under $(\zeta,\lambda)\mapsto\arg\bar\zeta\lambda$, as in \eqref{eq:necessary}; by \eqref{eq:necessary} and \eqref{eq:fourier}, $\psi_\omega(n)=2\int(1-\cos n\theta)\,d\tilde\nu$ for $n\ge1$. For $\theta\in\TT\setminus\{0\}$, $\frac1N\sum_{n=1}^N(1-\cos n\theta)\to1$, because $|\sum_{n=1}^N\cos n\theta|\le1/|\sin(\theta/2)|$. If $\psi_\omega(n)\to0$, the Ces\`aro means of $\psi_\omega$ tend to $0$, and Fatou's lemma gives $\tilde\nu=0$. Then $\psi_\omega(1)=2\|\omega\|_1=0$, so $\omega=0$. If $\omega$ vanishes on $\{|z|>\rho\}$, then $\psi_\omega(n)\le2n^2\rho^{2n-2}\|\omega\|_1\to0$, so $\omega=0$.
	\end{proof}
	
	\subsection{Radial weights and negative definite functions}\label{subsec:negdef}
	
	A function $\psi:\ZZ\to\RR$ is \emph{negative definite} if $\psi(0)=0$, $\psi(-n)=\psi(n)$ and $\sum_{n,k}c_n\bar c_k\psi(n-k)\le0$ whenever $\sum c_n=0$; equivalently, $e^{-t\psi}$ is positive definite for every $t>0$ (Schoenberg). By the L\'evy--Khintchine formula on the circle \cite[Ch.~I]{bergforst1975}, such a function has the form $\psi(n)=\sigma n^2+2\int_{\TT\setminus\{0\}}(1-\cos n\theta)\,\tilde\nu(d\theta)$, with $\sigma\ge0$ and a symmetric measure $\tilde\nu\ge0$ such that $\int(1-\cos\theta)\,d\tilde\nu<\infty$; it is the characteristic exponent of a symmetric L\'evy process on $\TT$ \cite{bertoin1996levy}.
	
	The radial criterion can also be written in terms of the sequence $\psi_\omega$.
	
	\begin{proposition}\label{prop:negdef}
		Let $\omega\ge0$ be a radial weight with $\int_\DD\omega\,dA<\infty$. Then $\omega$ is a Douglas weight if and only if $\psi_\omega$ is a negative definite function on $\ZZ$.
	\end{proposition}
	\begin{proof}
		If $\omega$ is a Douglas weight, Theorem~\ref{thm:A} and \eqref{eq:fourier} give $\psi_\omega(n)=2\int(1-\cos n\theta)\,d\tilde\nu(\theta)$ with $\tilde\nu=\frac1\pi j_\omega\,d\theta\ge0$. For finitely many numbers $c_n$ with $\sum c_n=0$,
		$$\sum_{n,k}c_n\bar c_k\big(1-\cos(n-k)\theta\big)=-\frac12\Big|\sum_nc_ne^{in\theta}\Big|^2-\frac12\Big|\sum_nc_ne^{-in\theta}\Big|^2\le0,$$
		so each function $n\mapsto1-\cos n\theta$ is negative definite. Integrating against $\tilde\nu$, we see that $\psi_\omega$ is negative definite.
		
		Conversely, let $\psi_\omega$ be negative definite. It is real and even, and $\psi_\omega(0)=0$, so the L\'evy--Khintchine formula gives $\psi_\omega(n)=\sigma n^2+2\int_{\TT\setminus\{0\}}(1-\cos n\theta)\,\tilde\nu(d\theta)$ as above. Since $1-\cos n\theta\le n^2(1-\cos\theta)$ and $(1-\cos n\theta)/n^2\to0$, dominated convergence gives $\psi_\omega(n)/n^2\to\sigma$. On the other hand $\psi_\omega(n)/n^2=2\int|z|^{2|n|-2}\omega\,dA\to0$ by \eqref{eq:psi}. Hence $\sigma=0$, so \eqref{eq:psinu} holds with a symmetric positive measure $\tilde\nu$, and the argument leading to \eqref{eq:QJradial} shows that $\omega\in\WLtwo$; in particular $\omega$ is a Douglas weight.
	\end{proof}
	
	Among radial weights, the constant weight has an $L^\infty$-neighbourhood consisting of Douglas weights.
		\begin{theorem}\label{thm:ball}
		%\revbegin
		Let $\varphi$ be a bounded, radial, measurable function on $\DD$ with $\|\varphi\|_\infty\le\frac45$. Then
		$$\omega=1+\varphi\in\WLtwo\qquad\text{and}\qquad j_\omega(\theta)\ \ge\ \Big(1-\tfrac54\|\varphi\|_\infty\Big)\,\frac1{4\sin^2(\theta/2)}\ \ge0\quad(0<\theta\le\pi).$$
		%\revend
	\end{theorem}
	To prove Theorem~\ref{thm:ball}, we first establish the following result.
	\begin{lemma}\label{lem:calJ}
		%\revbegin
		For $0<\theta\le\pi$,
		$$\int_0^1\big|\partial_\theta^2P(r^2,\theta)\big|\,\frac{dr}r=j_1(\theta)+\frac{\cos^2\theta}{4\sin^2\theta}\,\mathbf 1_{\{\theta<\pi/2\}}.$$
		In particular
		$$j_1(\theta)\le\int_0^1\big|\partial_\theta^2P(r^2,\theta)\big|\,\frac{dr}r\le\frac54\,j_1(\theta).$$
		%\revend
	\end{lemma}
	\begin{proof}
		%\revbegin
		Put $x=r^2$, so that $dr/r=dx/(2x)$. By \eqref{eq:d2Pseries}, $\partial_\theta^2P(x,\theta)=-2\sum_{k\ge1}k^2x^k\cos k\theta$, hence for $0<s<1$
		\begin{align*}F(s):=\int_0^s\partial_\theta^2P(x,\theta)\,\frac{dx}x&=-2\sum_{k\ge1}k\,s^k\cos k\theta=-2\,\Re\frac q{(1-q)^2}\\&=-\frac{2\big(s(1+s^2)\cos\theta-2s^2\big)}{\mathsf D(s,\theta)^2},\qquad q=se^{i\theta},\end{align*}
		because $q(1-\bar q)^2=q-2|q|^2+|q|^2\bar q$ and $|1-q|^2=\mathsf D(s,\theta)$. By \eqref{eq:Ctheta}, $|\partial_\theta^2P(x,\theta)|/x$ is bounded on $(0,1)$, so $F(s)\to F(1)=\frac1{1-\cos\theta}=2j_1(\theta)$ as $s\uparrow1$.
		By \eqref{eq:d2Pseries}, $\partial_\theta^2P=2x(1-x^2)\mathcal N(x)\mathsf D^{-3}$ with $\mathcal N(x)=4x\sin^2\theta-\mathsf D\cos\theta=-x^2\cos\theta+2(1+\sin^2\theta)\,x-\cos\theta$. If $\cos\theta\le0$, then $\mathcal N>0$ on $(0,1)$, so $\partial_\theta^2P\ge0$ and the integral equals $\frac12F(1)=j_1(\theta)$. If $\cos\theta>0$, then $\mathcal N(0)<0<\mathcal N(1)$ and the product of the two roots of $\mathcal N$ is $1$, so $\mathcal N$ has exactly one zero $u\in(0,1)$, with $\mathcal N<0$ on $(0,u)$ and $\mathcal N>0$ on $(u,1)$. Hence
		$$\int_0^1|\partial_\theta^2P|\,\frac{dx}x=F(1)-2F(u).$$
		From $\mathcal N(u)=0$ we get $u+u^{-1}=2(1+\sin^2\theta)/\cos\theta$. Therefore $\mathsf D(u,\theta)=u\,(u+u^{-1}-2\cos\theta)=4u\sin^2\theta/\cos\theta$ and $u(1+u^2)\cos\theta-2u^2=u^2\big((u+u^{-1})\cos\theta-2\big)=2u^2\sin^2\theta$, which gives $F(u)=-\cos^2\theta/(4\sin^2\theta)$. The integral in the statement is $\frac12\big(F(1)-2F(u)\big)=j_1(\theta)+\frac{\cos^2\theta}{4\sin^2\theta}$. Finally, for $0<\theta<\frac\pi2$,
		$$\frac{\cos^2\theta}{4\sin^2\theta}\Big/j_1(\theta)=\frac{\cos^2\theta}{4\cos^2(\theta/2)}\le\frac14,\qquad\text{since }\cos\theta\le\cos(\theta/2).$$
		%\revend
	\end{proof}
	\begin{proof}[Proof of Theorem \ref{thm:ball}]
		%\revbegin
		By linearity of \eqref{eq:jomega} in $\tilde\omega$, $j_\omega=j_1+j_\varphi$, and by Lemma~\ref{lem:calJ}, $|j_\varphi(\theta)|\le\|\varphi\|_\infty\int_0^1|\partial_\theta^2P(r^2,\theta)|\,\frac{dr}r\le\frac54\|\varphi\|_\infty\,j_1(\theta)$. Hence $j_\omega\ge(1-\frac54\|\varphi\|_\infty)\,j_1\ge0$, and Theorem~\ref{thm:A} applies.
	\end{proof}
	
	\subsection{Weights with a bounded Laplacian}\label{subsec:smooth}
	
	When a weight is bounded and has a bounded Laplacian, Theorem~\ref{thm:B} can be written as a pointwise inequality.
	
	\begin{proposition}\label{prop:smooth}
		Let $\omega\in C^2(\DD)$ with $0\le\omega\le C$ and $\Delta\omega$ bounded. Then the radial limit $h(\zeta)=\lim_{r\to1}\omega(r\zeta)$ exists for $m$-almost every $\zeta\in\TT$, and $\omega$ is a Douglas weight if and only if
		\begin{equation}\label{eq:pointwise}
			\int_\DD\Delta\omega(z)\,P_\zeta(z)P_\lambda(z)\,dA(z)\le\frac{4\big(h(\zeta)+h(\lambda)\big)}{|\zeta-\lambda|^2}\qquad\text{for $m\otimes m$-almost every }(\zeta,\lambda).
		\end{equation}
		In that case $\J_\omega=\mathfrak j\,dm\,dm$, where $4\mathfrak j(\zeta,\lambda)$ is the right-hand side minus the left-hand side of \eqref{eq:pointwise}.
	\end{proposition}
	\begin{proof}
		Let $\mu=-\frac14\Delta\omega\,dA=-\frac1{4\pi}\Delta\omega\,dx\,dy$, a signed measure with bounded density. Since $G[dA]=1-|z|^2$, we have $|G[\mu]|\le\frac14\|\Delta\omega\|_\infty(1-|z|^2)$, and $\Delta G[\mu]=\Delta\omega$ in $\mathcal D'(\DD)$. So $\omega-G[\mu]$ is harmonic by Weyl's lemma and bounded, hence equal to $P[h]$ for some $h\in L^\infty(\TT)$, and $h$ is the radial limit of $\omega$ almost everywhere because $G[\mu]\to0$ uniformly at $\TT$. Write $\mu=\mu^+-\mu^-$ and $h=h^+-h^-$. Then $\omega=\omega_1-\omega_2$ with $\omega_1=G[\mu^+]+P[h^+]$ and $\omega_2=G[\mu^-]+P[h^-]$, two bounded positive superharmonic weights. By Proposition~\ref{prop:superharmonic}, with $A_{\mu^\pm}$ defined by \eqref{eq:Amu}, and the formula $\J=-\frac12K$ of Theorem~\ref{thm:main}, $K_{\omega_1}$ and $K_{\omega_2}$ are, on $\TT^2\setminus\diag$, the functions $-2\big(A_{\mu^\pm}+(h^\pm(\zeta)+h^\pm(\lambda))|\zeta-\lambda|^{-2}\big)$. As $K_\omega=K_{\omega_1}-K_{\omega_2}$, we get $K_\omega=-2\mathfrak j$ on $\TT^2\setminus\diag$, with
		$$\mathfrak j(\zeta,\lambda)=-\frac14\int_\DD\Delta\omega\,P_\zeta P_\lambda\,dA+\frac{h(\zeta)+h(\lambda)}{|\zeta-\lambda|^2},$$
		a locally integrable function on $\TT^2\setminus\diag$ by Example~\ref{lem:AdA}. By Theorem~\ref{thm:main} and Lemma~\ref{lem:elem}(a), $\omega$ is a Douglas weight if and only if $K_\omega\le0$ on $\TT^2\setminus\diag$, that is, $\mathfrak j\ge0$ almost everywhere; then $\J_\omega=-\frac12K_\omega=\mathfrak j\,dm\,dm$.
	\end{proof}
	
	\begin{corollary}\label{cor:explicit}
		Let $\omega$ be as in Proposition~\ref{prop:smooth}, and assume that $h\ge h_0>0$ almost everywhere. If $\|\Delta\omega\|_\infty\le4h_0$, then $\omega$ is a Douglas weight.
	\end{corollary}
	\begin{proof}
		By Example~\ref{lem:AdA}(b),(d), $\big|\int_\DD\Delta\omega\,P_\zeta P_\lambda\,dA\big|\le\|\Delta\omega\|_\infty\mathsf A(\arg\bar\zeta\lambda)\le2\|\Delta\omega\|_\infty|\zeta-\lambda|^{-2}$, while the right-hand side of \eqref{eq:pointwise} is at least $8h_0|\zeta-\lambda|^{-2}$.
	\end{proof}
	
	With the numerical value $2\log2-1$ of $\sup_\theta\sin^2(\theta/2)\mathsf A(\theta)$, the same argument would allow $\|\Delta\omega\|_\infty\le2h_0/(2\log2-1)\approx5.18\,h_0$.
	\section{Examples and counterexamples}\label{sec:examples}
	
	We apply Theorems~\ref{thm:A} and~\ref{thm:B}, together with
	the results of Section~\ref{sec:further}, to several specific
	weights. Recall that a weight is called superharmonic if
	$\Delta\omega\leq 0$ in $\mathcal D'(\DD)$, and subharmonic
	if $\Delta\omega\geq 0$. 
	
	\subsection{Radial weights}\label{subsec:radialex}
	
	Powers of $|z|$, as well as their positive linear combinations,
	provide examples of radial subharmonic Douglas weights.
	
	\begin{example}\label{prop:subharm}
		\begin{enumerate}
			\item[(a)] For every $\gamma\ge0$, $|z|^{2\gamma}\in\WLtwo$. For $\gamma>0$ this weight is subharmonic and not superharmonic.
			\item[(b)] Let $\rho\ge0$ be a measure on $[0,\infty)$ with $\int\rho(d\gamma)/(\gamma+1)<\infty$, and $\omega(z)=\int_{[0,\infty)}|z|^{2\gamma}\rho(d\gamma)$. Then $\omega\in\WLtwo$, and $\omega$ is subharmonic and not superharmonic unless $\rho$ is carried by $\{0\}$.
			\item[(c)] For $0<\alpha<1$, the weights $(1-|z|^2)^{-\alpha}=\sum_{k\ge0}\frac{(\alpha)_k}{k!}|z|^{2k}$ and $-\log(1-|z|^2)=\sum_{k\ge1}\frac{|z|^{2k}}k$ are Douglas weights, subharmonic and not superharmonic.
		\end{enumerate}
	\end{example}
	
	\begin{example}\label{prop:power}
		The weight $(1-|z|^2)^s$, $s>-1$, is a Douglas weight if and only if $-1<s\le1$.
	\end{example}
	
	Small radial perturbations of the constant weight are Douglas weights by Theorem~\ref{thm:ball}, but need not be either superharmonic or subharmonic.
	
	\begin{example}\label{ex:wild}
		Let $N>\pi$ and $0<|\varepsilon|\le\frac4{5}$. By Theorem~\ref{thm:ball}, $\omega=1+\varepsilon\cos(N|z|^2)$ is a Douglas weight. For $f(z)=\cos(N|z|^2)$ and $r=|z|$, $\Delta f=-4N\sin(Nr^2)-4N^2r^2\cos(Nr^2)$, which equals $4\pi N>0$ on the circle $r^2=\pi/N$ and is negative for small $r>0$. Hence $\Delta\omega=\varepsilon\Delta f$ takes both signs.
	\end{example}
	
	\subsection{Non-radial weights}\label{subsec:nonradial}
	
	The simplest non-radial examples are obtained by addition.
	Let $S$ be a positive superharmonic weight and let $R$ be
	one of the weights considered in Example~\ref{prop:subharm}.
	Since $\WDou$ is a convex cone, $S+R$ is again a Douglas
	weight. For instance, $P_{\zeta_0}+|z|^2$ is a non-radial
	Douglas weight with Laplacian equal to $4$.
	
	For bounded weights with bounded Laplacian, Proposition~\ref{prop:smooth}
	and Corollary~\ref{cor:explicit} provide an explicit sufficient
	condition for membership in $\WDou$.
	
	\begin{example}\label{ex:r3cos}
		Let $\omega(z)=1+\varepsilon\,\Re z\,|z|^2=1+\varepsilon r^3\cos t$ with $0<\varepsilon\le\frac13$. Then $\Delta\omega=8\varepsilon\,\Re z$, so $\omega$ is neither superharmonic nor subharmonic, and $h(e^{i\alpha})=1+\varepsilon\cos\alpha\ge1-\varepsilon$. Since $\|\Delta\omega\|_\infty=8\varepsilon\le4(1-\varepsilon)$, Corollary~\ref{cor:explicit} shows that $\omega$ is a Douglas weight.
	\end{example}
	
	For weights with only the angular modes $0$ and $\pm1$, the matrix $M_\omega$ has three nonzero diagonals, and $K_\omega$ can be computed exactly off the diagonal.
	
	\begin{example}\label{lem:twomodes}
		Let $\omega$ be a weight with $\omega_j=0$ for $|j|\ge2$, so that $\omega(re^{it})=\tilde\omega(r)+2\Re\big(\omega_1(r)e^{it}\big)$. For $\theta\in\TT\setminus\{0\}$ put $c(\theta)=\min_{0\le x\le1}|1-xe^{i\theta}|>0$ and
		$$\kappa(\theta)=2e^{i\theta/2}\int_0^1\frac{\omega_1(r)}{r}\,L(r,\theta)\,dr,\qquad L(r,\theta)=4\sum_{m\ge1}m(m+1)\,r^{2m+1}\cos\big((m+\tfrac12)\theta\big)=8r\,\Re\frac{e^{i\theta/2}q}{(1-q)^3},$$
		where $q=r^2e^{i\theta}$. Then $|L(r,\theta)|\le8c(\theta)^{-3}r^3$, $\kappa$ is continuous on $\TT\setminus\{0\}$, and on $\TT^2\setminus\diag$ the distribution $K_\omega$ is the continuous function
		\begin{equation}\label{eq:Ktwo}
			K_\omega\big(e^{i\alpha},e^{i(\alpha+\theta)}\big)=-4j_\omega(\theta)+2\Re\big(\kappa(\theta)e^{i\alpha}\big).
		\end{equation}
		Consequently $\omega$ is a Douglas weight if and only if $|\kappa(\theta)|\le2j_\omega(\theta)$ for every $\theta\in(0,\pi]$.
	\end{example}
	
	\begin{example}\label{ex:cos}
		For $0<\varepsilon\le1$ let $\sigma_\varepsilon(z)=1+\varepsilon\cos(\arg z)=1+\varepsilon\,\Re z/|z|$. Then $\sigma_\varepsilon\ge0$, its radialization is the constant $1$, and $\omega_1=\varepsilon/2$, so Example~\ref{lem:twomodes} applies with $j_{\sigma_\varepsilon}=j_1$ and $\kappa(\theta)=\varepsilon e^{i\theta/2}I(\theta)$, where $I(\theta)=\int_0^1L(r,\theta)\,\frac{dr}r$. For $0<s<1$, termwise integration and $4m(m+1)=(2m+1)^2-1$ give
		$$\int_0^s L(r,\theta)\,\frac{dr}r=\sum_{m\ge1}\Big((2m+1)-\frac1{2m+1}\Big)s^{2m+1}\cos\big((m+\tfrac12)\theta\big)=\Re\Big[\frac{w(1+w^2)}{(1-w^2)^2}-\frac12\log\frac{1+w}{1-w}\Big],$$
		with $w=se^{i\theta/2}$; the terms with $m=0$ of the two series are equal, so they cancel. Let $s\uparrow1$. The left-hand side tends to $I(\theta)$ by dominated convergence, since $|L(r,\theta)|/r\le8c(\theta)^{-3}r^2$. On the right, for $w=e^{i\theta/2}$ and $0<\theta\le\pi$, $\frac{w(1+w^2)}{(1-w^2)^2}=-\frac{\cos(\theta/2)}{2\sin^2(\theta/2)}$ and $\frac{1+w}{1-w}=i\cot\frac\theta4$. Hence
		$$I(\theta)=-a(\theta),\qquad a(\theta)=\frac{\cos(\theta/2)}{2\sin^2(\theta/2)}+\frac12\log\cot\frac\theta4\ge0,$$
		and the condition $|\kappa|\le2j_1=\frac1{2\sin^2(\theta/2)}$ of Example~\ref{lem:twomodes} becomes
		\begin{equation}\label{eq:thr}
			\varepsilon\,g(\theta)\le1\quad\text{for all }\theta\in(0,\pi],\qquad g(\theta)=2\sin^2\tfrac\theta2\,a(\theta)=\cos\frac\theta2+\sin^2\frac\theta2\,\log\cot\frac\theta4 .
		\end{equation}
		The function $g$ is continuous on $(0,\pi]$, $g(\pi)=0$, and $g(\theta)=1+\frac14\theta^2\big(\log\frac4\theta-\frac12\big)+o(\theta^2)$ as $\theta\to0$. Hence $\sup g>1$, the supremum is attained, and $\sigma_\varepsilon$ is a Douglas weight exactly when $\varepsilon\le\varepsilon^*=1/\sup g<1$. For $\varepsilon^*<\varepsilon\le1$, $\sigma_\varepsilon$ is not a Douglas weight, although its radialization is. Numerically, $\sup g\approx1.1997$ is attained at $\theta\approx1.1706$, and $\varepsilon^*\approx0.8336$.
	\end{example}
	
	\subsection{Counterexamples}\label{subsec:outside}
	
	Corollary~\ref{cor:boundary} and Example~\ref{ex:cos} already give positive weights that are not Douglas weights: $(1-|z|^2)^s$ with $s>1$, the nonzero weights with compact support in $\DD$, and $\sigma_\varepsilon$ for $\varepsilon^*<\varepsilon\le1$. A large constant turns a compactly supported weight into a Douglas weight.
	
	\begin{proposition}\label{prop:nonsolid}
		Let $\omega_0\not\equiv0$ be a weight that vanishes almost everywhere on $\{|z|>\rho\}$ for some $\rho<1$. Then $\omega_0\notin\WDou$, whereas $c+\omega_0\in\WDou$ for every constant $c\ge4(1-\rho)^{-4}\|\omega_0\|_1$. If moreover $\omega_0$ is smooth, then $c+\omega_0$ is not superharmonic.
	\end{proposition}
	\begin{proof}
		The first assertion is Corollary~\ref{cor:boundary}. For $|z|\le\rho$, \eqref{eq:kappaz} gives $|\kappa_z(\zeta,\lambda)|\le4|\zeta-z|^{-2}|\lambda-z|^{-2}\le4(1-\rho)^{-4}$, so $K_{\omega_0}=\int\omega_0\,\kappa_z\,dA$ is a continuous function on $\TT^2$ bounded by $4(1-\rho)^{-4}\|\omega_0\|_1$. By Example~\ref{ex:j1}, $K_c=-2\J_c=-4c\,|\zeta-\lambda|^{-2}\le-c$ off the diagonal. Hence $K_{c+\omega_0}\le0$ off the diagonal, and $c+\omega_0\in\WDou$ by Theorem~\ref{thm:main}. If $\omega_0$ is smooth, it has compact support, so $\int_\DD\Delta\omega_0\,dA=0$; if $\Delta\omega_0\le0$, continuity would give $\Delta\omega_0\equiv0$, and $\omega_0$ would be a compactly supported harmonic function, hence $0$. So $\Delta(c+\omega_0)=\Delta\omega_0$ is positive somewhere.
	\end{proof}
	
	The stability of the constant weight under small radial perturbations (Theorem~\ref{thm:ball}) has no non-radial analogue.
	
	\begin{proposition}\label{prop:nostab}
		For every $\varepsilon\in(0,1]$ there is a measurable function $\phi$ on $\DD$ with $|\phi|\le1$ such that the weight $1+\varepsilon\phi$ is not a Douglas weight.
	\end{proposition}
	\begin{proof}
		Fix $\zeta_0\ne\lambda_0$ in $\TT$ and put $d=|\zeta_0-\lambda_0|$. We first show that
		\begin{equation}\label{eq:kappadiv}
			\int_\DD|\kappa_z(\zeta_0,\lambda_0)|\,dA(z)=+\infty .
		\end{equation}
		With the notation of the proof of Lemma~\ref{lem:elem}, $\kappa_z(\zeta_0,\lambda_0)=\Re\big(F_{\zeta_0}'(z)\overline{F_{\lambda_0}'(z)}\big)$ with $F_\zeta'(z)=2\zeta/(\zeta-z)^2$. Write $z=\zeta_0(1-\varrho e^{i\eta})$ with $\varrho>0$ and $|\eta|<\frac\pi2$. Then $F_{\zeta_0}'(z)=2\bar\zeta_0\,e^{-2i\eta}\varrho^{-2}$, and $F_{\lambda_0}'(z)=-2\bar\zeta_0\,d^{-2}+O(\varrho)$ as $\varrho\to0$, because $(\lambda_0-\zeta_0)^2=-\zeta_0\lambda_0\,d^2$. Hence
		$$\kappa_z(\zeta_0,\lambda_0)=-\frac{4\cos2\eta}{d^2\varrho^2}+O\Big(\frac1\varrho\Big)\qquad(\varrho\to0).$$
		On the sector $\{|\eta|<\frac\pi8,\ 0<\varrho<\varrho_0\}$, contained in $\DD$ for small $\varrho_0$, this gives $|\kappa_z(\zeta_0,\lambda_0)|\ge c_0\varrho^{-2}$ with $c_0>0$, and \eqref{eq:kappadiv} follows because $dA=\frac1\pi\varrho\,d\varrho\,d\eta$ there.
		
		Let $\varphi_k\ge0$ be functions in $C_c^\infty(\TT^2\setminus\diag)$ with $\iint\varphi_k\,dm\,dm=1$ whose supports shrink to $(\zeta_0,\lambda_0)$, and put $\Phi_k(z)=\langle\kappa_z,\varphi_k\rangle=\frac12\Delta(\mathcal P\varphi_k)(z)$, a bounded continuous function on $\DD$ (Lemma~\ref{lem:elem}). Since $\kappa_z$ is continuous in $(\zeta,\lambda)$, $\Phi_k(z)\to\kappa_z(\zeta_0,\lambda_0)$, and Fatou's lemma with \eqref{eq:kappadiv} gives $\int_\DD|\Phi_k|\,dA\to+\infty$. On the other hand, by Lemma~\ref{lem:elem} and Example~\ref{ex:j1},
		$$\int_\DD\Phi_k\,dA=\langle K_1,\varphi_k\rangle=-4\iint\frac{\varphi_k(\zeta,\lambda)}{|\zeta-\lambda|^2}\,dm\,dm\longrightarrow-\frac4{d^2}.$$
		Choose $k$ with $\varepsilon\int_\DD|\Phi_k|\,dA>4\iint\varphi_k|\zeta-\lambda|^{-2}dm\,dm$, and let $\phi$ be the sign of $\Phi_k$. Then $\omega=1+\varepsilon\phi$ satisfies $1-\varepsilon\le\omega\le1+\varepsilon$ and
		$$\frac12\int_\DD\omega\,\Delta(\mathcal P\varphi_k)\,dA=\int_\DD\Phi_k\,dA+\varepsilon\int_\DD|\Phi_k|\,dA>0 ,$$
		so $\omega$ does not satisfy condition~\eqref{eq:weaksuper}. By Theorem~\ref{thm:B} it is not a Douglas weight.
	\end{proof}
	
	The function $\phi$ used in the proof takes the values $\pm1$
	on sectors with vertices at $\zeta_0$ and $\lambda_0$.
	In contrast, for weights with bounded Laplacian,
	Corollary~\ref{cor:explicit} provides a sufficient smallness
	condition, as illustrated by Example~\ref{ex:r3cos}.
	
	Finally, a singular boundary measure provides an example
	distinguishing the two classes $\WLip$ and $\WLtwo$
	introduced in the introduction.
	
	\begin{proposition}\label{prop:strict}
		Let $\omega=G[\mu]+P[\nu]$ be a positive superharmonic weight whose boundary measure $\nu$ has a nonzero singular part $\nu_s$ with respect to $m$. Then $\omega\in\WLip\setminus\WLtwo$. In particular, for the Poisson kernel $\omega=P_{\zeta_0}$, $\J_\omega=(\delta_{\zeta_0}\otimes m+m\otimes\delta_{\zeta_0})|\zeta-\lambda|^{-2}$, and the function $f=\mathbf 1_{\{\zeta_0\}}$, which is $0$ in $L^2(\TT)$, satisfies $Q_\omega(f)=0$ and $Q_{\J_\omega}(f)=+\infty$.
	\end{proposition}
	\begin{proof}
		By Proposition~\ref{prop:superharmonic}, $\omega\in\WLip$ with $\J_\omega\ge(\nu\otimes m)|\zeta-\lambda|^{-2}$. Suppose that $\omega\in\WLtwo$ with a measure $\J'$. The identity $Q_\omega=Q_{\J'}$ holds in particular on Lipschitz functions, so $\J'=\J_\omega$ by the uniqueness proved at the end of Section~\ref{sec:radial}. Let $N$ be an $m$-null Borel set with $\nu_s(N)>0$. The Borel function $\mathbf 1_N$ is $0$ in $L^2(\TT)$, so $Q_\omega(\mathbf 1_N)=0$. However,
		$$Q_{\J_\omega}(\mathbf 1_N)\ge\J_\omega\big(N\times(\TT\setminus N)\big)\ge\int_N d\nu(\zeta)\int_\TT\frac{dm(\lambda)}{|\zeta-\lambda|^2}=+\infty,$$
		because $m(\TT\setminus N)=1$ and $\int_\TT|\zeta-\lambda|^{-2}dm(\lambda)=+\infty$ for every $\zeta$. This contradiction proves the first assertion; the second is the case $\mu=0$, $\nu=\delta_{\zeta_0}$, $N=\{\zeta_0\}$.
	\end{proof}

\bibliographystyle{abbrv}
	\bibliography{references}

@incollection{abakumov2015cyclicity,
  author    = {Abakumov, Evgueni and El-Fallah, Omar and Kellay, Karim and Ransford, Thomas},
  title     = {Cyclicity in the harmonic {D}irichlet space},
  booktitle = {Conference on Harmonic and Functional Analysis, Operator Theory and Applications},
  series    = {Theta Ser. Adv. Math.},
  volume    = {19},
  publisher = {Theta},
  address   = {Bucharest},
  year      = {2017},
  pages     = {1--10}
}

@phdthesis{aleman1993multiplication,
  author = {Aleman, Alexandru},
  title  = {The multiplication operator on {H}ilbert spaces of analytic functions},
  type   = {Habilitationsschrift},
  school = {FernUniversit{\"a}t in Hagen},
  year   = {1993}
}

@article{bahajji2020approximation,
  title={Approximation in spaces of analytic functions},
  author={Bahajji-El Idrissi, Hafid and El-Fallah, Omar},
  journal={Stud. Math.},
  volume={255},
  pages={209--217},
  year={2020}
}

@unpublished{bahajji2026superharmonically,
  author = {Bahajji-El Idrissi, H. and El-Fallah, O. and Elmadani, Y. and Hanine, A.},
  title  = {Superharmonically weighted {D}irichlet spaces},
  note   = {Preprint, arXiv:2605.13787},
  doi    = {10.48550/arXiv.2605.13787},
  year   = {2026}
}

@article{bao2018dirichlet,
  author  = {Bao, Guanlong and G{\"o}{\u{g}}{\"u}{\c{s}}, Nihat G{\"o}khan and Pouliasis, Stamatis},
  title   = {On {D}irichlet spaces with a class of superharmonic weights},
  journal = {Can. J. Math.},
  volume  = {70},
  number  = {4},
  pages   = {721--741},
  year    = {2018},
  doi     = {10.4153/CJM-2017-005-1}
}

@unpublished{BEHJL2026,
  author = {Bourabiaa, J. and Elmadani, Y. and Hanine, A. and Jamal, A. and Labghail, I.},
  title  = {Douglas weighted {D}irichlet spaces: analytic and probabilistic aspects},
  note   = {Preprint, arXiv:2607.26870},
  doi    = {10.48550/arXiv.2607.26870},
  year   = {2026}
}

@book{bergforst1975,
  author    = {Berg, Christian and Forst, Gunnar},
  title     = {Potential Theory on Locally Compact Abelian Groups},
  series    = {Ergebnisse der Mathematik und ihrer Grenzgebiete},
  volume    = {87},
  publisher = {Springer-Verlag},
  address   = {New York-Heidelberg},
  year      = {1975}
}

@book{bertoin1996levy,
  author    = {Bertoin, Jean},
  title     = {L{\'e}vy processes},
  series    = {Cambridge Tracts in Mathematics},
  volume    = {121},
  publisher = {Cambridge University Press},
  address   = {Cambridge},
  year      = {1996},
  pages     = {x + 265},
  isbn      = {978-0-521-56243-0},
  isbn13    = {978-0521562430},
  mrnumber  = {1406564},
  doi       = {10.1017/CBO9780511809781}
}

@article{beurling1958espaces,
  author  = {Beurling, A. and Deny, J.},
  title   = {Espaces de {D}irichlet. {I}. {L}e cas {\'e}l{\'e}mentaire},
  journal = {Acta Math.},
  volume  = {99},
  number  = {1},
  pages   = {203--224},
  year    = {1958}
}

@article{carleson1960representation,
  author  = {Carleson, Lennart},
  title   = {A representation formula for the {D}irichlet integral},
  journal = {Math. Z.},
  volume  = {73},
  pages   = {190--196},
  year    = {1960}
}

@article{chen2006traces,
  author  = {Chen, Zhen-Qing and Fukushima, Masatoshi and Ying, Jiangang},
  title   = {Traces of symmetric {M}arkov processes and their characterizations},
  journal = {Ann. Probab.},
  volume  = {34},
  number  = {3},
  pages   = {1052--1102},
  year    = {2006},
  doi     = {10.1214/009117905000000657}
}

@book{chenfukushima2012,
  author    = {Chen, Zhen-Qing and Fukushima, Masatoshi},
  title     = {Symmetric {M}arkov Processes, Time Change, and Boundary Theory},
  series    = {London Mathematical Society Monographs Series},
  volume    = {35},
  publisher = {Princeton University Press},
  address   = {Princeton, NJ},
  year      = {2012}
}

@article{ChevrotGuillotRansford2010,
  author  = {Nicolas Chevrot and Dominique Guillot and Thomas Ransford},
  title   = {De {B}ranges--{R}ovnyak spaces and {D}irichlet spaces},
  journal = {J. Funct. Anal.},
  volume  = {259},
  number  = {9},
  pages   = {2366--2383},
  year    = {2010},
  doi     = {10.1016/j.jfa.2010.06.018}
}

@article{douglas1931solution,
  author  = {Douglas, Jesse},
  title   = {Solution of the problem of {P}lateau},
  journal = {Trans. Am. Math. Soc.},
  volume  = {33},
  number  = {1},
  pages   = {263--321},
  year    = {1931}
}

@article{el2009brown,
  title={On the {B}rown--{S}hields conjecture for cyclicity in the {D}irichlet space},
  author={El-Fallah, Omar and Kellay, Karim and Ransford, Thomas},
  journal={Adv. Math.},
  volume={222},
  number={6},
  pages={2196--2214},
  year={2009}
}

@book{el2014primer,
  author    = {El-Fallah, Omar and Kellay, Karim and Mashreghi, Javad and Ransford, Thomas},
  title     = {A primer on the {D}irichlet space},
  series    = {Cambridge Tracts in Mathematics},
  volume    = {203},
  publisher = {Cambridge University Press},
  address   = {Cambridge},
  year      = {2014},
  isbn      = {978-1-107-04752-5},
  doi       = {10.1017/CBO9781107239425},
  lccn      = {2013372391},
  pages      = {xiii + 211},
  mrnumber   = {3184178}
}

@article{el2016cyclicity,
  title={Cyclicity and invariant subspaces in {D}irichlet spaces},
  author={El-Fallah, Omar and Elmadani, Youssef and Kellay, Karim},
  journal={J. Funct. Anal.},
  volume={270},
  number={9},
  pages={3262--3279},
  year={2016}
}

@article{el2016dirichlet,
  title={{D}irichlet spaces with superharmonic weights and {D}e {B}ranges--{R}ovnyak spaces},
  author={El-Fallah, Omar and Kellay, Karim and Klaja, Hubert and Mashreghi, Javad and Ransford, Thomas},
  journal={Complex Anal. Oper. Theory},
  volume={10},
  number={1},
  pages={97--107},
  year={2016}
}

@article{el2019kernel,
  title={Kernel and capacity estimates in {D}irichlet spaces},
  author={El-Fallah, O. and Elmadani, Y. and Kellay, K.},
  journal={J. Funct. Anal.},
  volume={276},
  number={3},
  pages={867--895},
  year={2019}
}

@article{el2022extremal,
  title={Extremal functions and invariant subspaces in {D}irichlet spaces},
  author={El-Fallah, O. and Elmadani, Y. and Labghail, I.},
  journal={Adv. Math.},
  volume={408},
  pages={108604},
  year={2022}
}

@article{el2024ahlfors,
  author  = {El-Fallah, O. and Elmadani, Y. and Labghail, I.},
  title   = {Ahlfors--{D}avid regular sets, point spectrum and {D}irichlet spaces},
  journal = {Results Math.},
  volume  = {79},
  number  = {2},
  pages   = {74},
  year    = {2024}
}

@article{elmadani2019cyclicity,
  title={Cyclicity in {D}irichlet spaces},
  author={Elmadani, Youssef and Labghail, Imane},
  journal={Can. Math. Bull.},
  volume={62},
  number={2},
  pages={247--257},
  year={2019}
}

@article{fukushima2004timechanges,
  author  = {Fukushima, Masatoshi and He, Ping and Ying, Jiangang},
  title   = {Time changes of symmetric diffusions and {F}eller measures},
  journal = {Ann. Probab.},
  volume  = {32},
  number  = {4},
  pages   = {3138--3166},
  year    = {2004},
  doi     = {10.1214/009117904000000649}
}

@book{fukushima2010dirichlet,
  author    = {Fukushima, Masatoshi and Oshima, Yoichi and Takeda, Masayoshi},
  title     = {{D}irichlet forms and symmetric {M}arkov processes},
  series    = {De Gruyter Studies in Mathematics},
  volume    = {19},
  edition   = {2},
  publisher = {Walter de Gruyter},
  address   = {Berlin/New York},
  year      = {2011},
  pages     = {viii + 489},
  isbn      = {978-3-11-021808-4},
  eisbn     = {978-3-11-021809-1},
  doi       = {10.1515/9783110218091},
  lccn      = {2010041939},
  mrnumber  = {2778606},
  zbl       = {1225.31001}
}

@article{gallardo2025zero,
  title={Zero-free regions of the {R}iemann zeta function and approximation in weighted {D}irichlet spaces},
  author={Gallardo-Guti{\'e}rrez, Eva A. and Seco, Daniel},
  journal={Complex Anal. Oper. Theory},
  volume={19},
  number={2},
  pages={38},
  year={2025}
}

@article{guillot2012fine,
  title={Fine boundary behavior and invariant subspaces of harmonically weighted {D}irichlet spaces},
  author={Guillot, Dominique},
  journal={Complex Anal. Oper. Theory},
  volume={6},
  number={6},
  pages={1211--1230},
  year={2012}
}

@book{hormander1983analysis,
  author = {H{\"o}rmander, Lars},
  title = {The analysis of linear partial differential operators. {I}},
  series = {Grundlehren der mathematischen Wissenschaften},
  volume = {256},
  publisher = {Springer-Verlag},
  address = {Berlin},
  year = {1983}
}

@article{idrissi2020blaschke,
  title={Blaschke products and zero sets in weighted {D}irichlet spaces},
  author={Bahajji-El Idrissi, Hafid and El-Fallah, Omar},
  journal={Potential Anal.},
  volume={53},
  number={4},
  pages={1299--1316},
  year={2020}
}

@article{idrissi2025douglas,
  title={Douglas-type formula for weighted {B}esov spaces},
  author={Bahajji-El Idrissi, Hafid and El-Fallah, Omar},
  journal={Collect. Math.},
  pages={1--17},
  year={2025}
}

@book{koosis1998introduction,
  author    = {Koosis, Paul},
  title     = {Introduction to {$H^p$} spaces},
  series    = {Cambridge Tracts in Mathematics},
  volume     = {115},
  edition    = {2},
  publisher  = {Cambridge University Press},
  address    = {Cambridge},
  year       = {1998},
  pages      = {xiv + 287},
  isbn       = {978-0-521-45521-3},
  doi        = {10.1017/CBO9780511470950},
  mrnumber   = {1669574}
}

@article{richter1991formula,
  author  = {Richter, Stefan and Sundberg, Carl},
  title   = {A formula for the local {D}irichlet integral},
  journal = {Mich. Math. J.},
  volume  = {38},
  number  = {3},
  pages   = {355--379},
  year    = {1991}
}

@article{richter1991representation,
  title={A representation theorem for cyclic analytic two-isometries},
  author={Richter, Stefan},
  journal={Trans. Am. Math. Soc.},
  volume={328},
  number={1},
  pages={325--349},
  year={1991}
}

@article{ross1994hyperinvariant,
  author  = {Richter, Stefan and Ross, William T. and Sundberg, Carl},
  title   = {Hyperinvariant subspaces of the harmonic {D}irichlet space},
  journal = {J. Reine Angew. Math.},
  volume  = {448},
  pages   = {1--26},
  year    = {1994}
}

@article{Sarason1997,
  author  = {Donald Sarason},
  title   = {Local {D}irichlet spaces as {D}e {B}ranges--{R}ovnyak spaces},
  journal = {Proc. Am. Math. Soc.},
  volume  = {125},
  number  = {7},
  pages   = {2133--2139},
  year    = {1997}
}

@article{shimorin2001reproducing,
  title={Reproducing kernels and extremal functions in {D}irichlet--type spaces},
  author={Shimorin, Serguei},
  journal={J. Math. Sci.},
  volume={107},
  number={4},
  pages={4108--4124},
  year={2001}
}

@article{shimorin2002complete,
  title={Complete {N}evanlinna--{P}ick property of {D}irichlet--type spaces},
  author={Shimorin, Serguei},
  journal={J. Funct. Anal.},
  volume={191},
  number={2},
  pages={276--296},
  year={2002}
}
	
\end{document}